\documentclass[12pt,reqno]{amsart}
\usepackage{amsmath,amsfonts,amsthm,amsopn,amssymb}
\usepackage{cite,marginnote}
\usepackage{color,enumitem,graphicx}
\usepackage[colorlinks=true,urlcolor=blue,
citecolor=red,linkcolor=blue,linktocpage,pdfpagelabels,
bookmarksnumbered,bookmarksopen]{hyperref}
\usepackage[english]{babel}

\usepackage[left=2.9cm,right=2.9cm,top=2.8cm,bottom=2.8cm]{geometry}

\numberwithin{equation}{section}

\makeindex

\newtheorem{Theorem}{Theorem}[section]

\newtheorem{Remark}{Remark}[section]

\newtheorem{Definition}{Definition}[section]

\newtheorem{lemma}{Lemma}[section]
\newtheorem{iteration lemma}{iteration Lemma}[section]

\newcommand{\s}{\section}

\newcommand{\R}{\mathbb R}

\newcommand{\bt}{\begin{theorem}}
\newcommand{\et}{\end{theorem}}
\newcommand{\bl}{\begin{lemma}}
\newcommand{\el}{\end{lemma}}
\newcommand{\bd}{\begin{definition}}
\newcommand{\ed}{\end{definition}}
\newcommand{\bc}{\begin{corollary}}
\newcommand{\ec}{\end{corollary}}
\newcommand{\bp}{\begin{proof}}
\newcommand{\ep}{\end{proof}}
\newcommand{\bx}{\begin{example}}
\newcommand{\ex}{\end{example}}
\newcommand{\bi}{\begin{exercise}}
\newcommand{\ei}{\end{exercise}}
\newcommand{\bo}{\begin{proposition}}
\newcommand{\eo}{\end{proposition}}
\newcommand{\br}{\begin{remark}}
\newcommand{\er}{\end{remark}}
\newcommand{\beq}{\begin{equation}}
\newcommand{\eeq}{\end{equation}}
\newcommand{\ba}{\begin{align}}
\newcommand{\ea}{\end{align}}
\newcommand{\bn}{\begin{enumerate}}
\newcommand{\en}{\end{enumerate}}
\newcommand{\bg}{\begin{align*}}
\newcommand{\bcs}{\begin{cases}}
\newcommand{\ecs}{\end{cases}}

\newcommand{\bean}{\begin{eqnarray*}}
\newcommand{\eean}{\end{eqnarray*}}

\def\R{\mathbb{R}}

\def\bd{\mathrm{bd}\,}

\title[Existence and multiplicity of sign-changing solutions]{Existence and multiplicity of sign-changing solutions for quasilinear Schr\"{o}dinger equations with sub-cubic nonlinearity}

\author[H. Zhang]{Hui Zhang}
\author[Z. S. Liu]{Zhisu Liu}
\author[C.-L. Tang]{Chun-Lei Tang}
\author[J. J. Zhang]{Jianjun Zhang}

\address[H.\ Zhang]{\newline\indent Department of  Mathematics,
Jinling Institute of Technology,
\newline\indent
Nanjing 211169, China
\newline\indent and
\newline\indent Department of Mathematics,
Nanjing University,
\newline\indent
Nanjing 210093, China}
\email{\href{mailto:huihz0517@126.com}{huihz0517@126.com}}

\address[Z. S.\ Liu]{\newline\indent Center for Mathematical Sciences, China University of Geosciences,
\newline\indent
Wuhan 430074, China}
\email{\href{mailto:liuzhisu@cug.edu.cn}{liuzhisu@cug.edu.cn}}

\address[C.-L.\ Tang]{\newline\indent School of Mathematics and Statistics, Southwest University,
\newline\indent
Chongqing 400715, China}
\email{\href{mailto:tangcl@swu.edu.cn}{tangcl@swu.edu.cn}}

\address[J. J. \ Zhang]{\newline\indent College of Mathematics and Statistics, Chongqing Jiaotong University,
\newline\indent
Chongqing 400074, China}
\email{\href{mailto:zhangjianjun09@tsinghua.org.cn}{zhangjianjun09@tsinghua.org.cn}}

\thanks{Corresponding author: Chun-Lei Tang ({\tt tangcl@swu.edu.cn})}

\thanks{Hui Zhang was supported by China Postdoctoral Science Foundation (No. 2021M691527). Zhisu Liu was supported by the NSFC (No.~11701267), the Hunan Natural Science Excellent Youth Fund (No.~2020JJ3029) and the Fundamental Research Funds for the Central
Universities, China University of Geosciences (Wuhan, Grant number: CUGST2). Chun-Lei Tang was supported by the NSFC (No. 11971393). Jianjun Zhang was supported by the NSFC (No. 11871123) and Team Building Project for Graduate Tutors in Chongqing(JDDSTD201802).}


\begin{document}

\begin{abstract}
In this paper, we consider the
 quasilinear Schr\"{o}dinger equation
\begin{equation*}
-\Delta u+V(x)u-u\Delta(u^2)=g(u),\ \
x\in \mathbb{R}^{3},
\end{equation*}
where $V$ and $g$ are continuous functions. Without the coercive condition on $V$ or the monotonicity condition on $g$, we show that the problem above has a least energy sign-changing solution and infinitely many sign-changing solutions.  Our results especially solve the problem above in the case where $g(u)=|u|^{p-2}u$ ($2<p<4$) and complete some recent related works on sign-changing solutions, in the sense that, in the literature only the case $g(u)=|u|^{p-2}u$ ($p\geq4$) was considered. The main results in the present paper are obtained by a new perturbation approach and the method of invariant sets of descending flow. In addition, in some cases where the functional merely satisfies the Cerami condition, a deformation lemma under the Cerami condition is developed.
\end{abstract}

\maketitle

\s{Introduction and main results}
\renewcommand{\theequation}{1.\arabic{equation}}
The present paper is devoted to studying the existence and multiplicity of sign-changing solutions of the quasilinear Schr\"{o}dinger equation
\begin{equation}\label{1.0.0}
-\triangle u+V(x)u-u \Delta (u^2)={g}(u), \ \  x\in \mathbb{R}^3,\hskip 3 cm
\end{equation}
where $V\in C(\mathbb{R}^3,\mathbb{R})$ and $g\in(\mathbb{R},\mathbb{R})$. Problem (\ref{1.0.0}) is referred to as the so-called Modified Nonlinear Schr\"odinger Equation
(MNLS) and is related to solitary wave solutions of the equation
\begin{equation}\label{1.1}
i\partial_tz=-\triangle z+W(x)z-m(|z|^2)z-\kappa\triangle
\phi(|z|^2)\phi'(|z|^2)z, \hskip 3 cm
\end{equation}
 where $z:\mathbb{R}\times \mathbb{R}^3\rightarrow \mathbb{C},$
$W:\mathbb{R}^3\rightarrow \mathbb{R}$ is a given potential,
$m,\phi:\mathbb{R}^+\rightarrow \mathbb{R}$ are suitable functions, and
$\kappa\in\R$. The form of (\ref{1.1}) has been derived as models of several physical phenomena
corresponding to various types of $\kappa$ and $\phi$. For example, if $\kappa=0,$ (\ref{1.1}) turns out to be a semilinear Schr\"{o}dinger equation, which has been widely
investigated, we refer the readers to \cite{BERE,Rap,WM1}. The case
$\phi(s)=s$, as a model of the time evolution of the condensate wave
function in super-fluid film, has been studied by Kurihara in
\cite{Ku}. While for $\phi(s)=\sqrt{1+s},$ the equations are the models
of the self-channeling of a high-power ultra short laser in matter, see \cite{ChenS}. For more physical
applications, we refer to \cite{BL,LPT,PS} and references therein.

In the past decades, great progress has been made in studying sign-changing solutions of (\ref{1.0.0}) and such problems have attracted a considerable researchers' attention. In \cite{LWW}, Liu et al.
considered the equation
\begin{equation}\label{1.3}
-\triangle u+V(x)u-u \Delta (u^2)=|u|^{p-2}u, \ \  x\in \mathbb{R}^N,
\end{equation}
where $N\geq3$, $p\in[4,22^*)$ with $2^*=\frac{2N}{N-2}$, and $V(x)\in C(\mathbb{R}^N)$ satisfies
\vskip 0.2 true cm
\noindent(V$'_1$)\ $0<\inf_{\mathbb{R}^N}V(x)\leq \lim_{|x|\rightarrow+\infty}V(x):=V_\infty$,\ \text{and} \ $V(x)\leq V_\infty-\frac{A}{1+|x|^m}$ \text{for} \ $|x|\geq M$,
\vskip 0.2 true cm
\noindent where $A, M, m$ are positive constants. They proved that (\ref{1.3}) has a least energy sign-changing solution by using an approximating sequence of problems on a Nehari manifold defined in an appropriate subset of $H^1(\mathbb{R}^N)$. In \cite{DPW1}, Deng et al. treated the equation (\ref{1.3}) with $p\in(4,22^*)$ and showed that, for any given $k\in\mathbb{N}$, there is a pair of
sign-changing solutions with $k$ nodes by using a minimization argument and  an energy comparison method. Later, Deng et al. \cite{DPW2} extended the results in \cite{DPW1} to the critical growth case. Zhang et al. \cite{ZL} showed the existence of infinitely many sign-changing solutions of the equation
 \begin{equation}\label{1.5}
-\triangle u+u-u \Delta (u^2)=a(x)|u|^{p-2}u, \ \  x\in \mathbb{R}^N,
\end{equation}
where $p\in(4,22^*)$ and $a(x)$ satisfies
$$a(x)>0,\ a\in L^r(\mathbb{R}^N)\ \text{with}\ r\geq{{22^*}/{(22^*-p)}},$$
and the proof is based on the methods of perturbation and invariant sets of descending
flow.  Recently, Yang et al. \cite{YANG1} dealt with the equation with critical or supercritical growth
$$-\triangle u+V(x)u-u \Delta (u^2)=a(x)[g(u)+|u|^{p-2}u], \ \  x\in \mathbb{R}^N,$$
where $p\geq 22^*$.  By assuming that $a(x)>0$ a.e. in $\mathbb{R}^N$, $V$ satisfies one of the following conditions
\vskip 0.2 true cm
\noindent(V$'_2$) $V(x)\geq V_0>0$ for all $x\in\mathbb{R}^N$,\ $V(x)=V(|x|)$ and $V\in L^\infty(\mathbb{R}^N)$,
\vskip 0.1 true cm
\noindent(V$'_3$) $V(x)\geq V_0>0$ for all $x\in\mathbb{R}^N$,\ $\lim_{|x|\rightarrow\infty}V(x)=+\infty$,
\vskip 0.1 true cm
\noindent and $g\in C(\mathbb{R})$ is odd, subcritical at infinity, superlinear near zero and satisfies
\vskip 0.1 true cm
\noindent($g'_1$) $\frac{g(t)}{t^3}$ \ {is nondecreasing  in} $t\neq0,$ \vskip 0.1 true cm
\noindent($g'_2$) $g(t)t\geq\mu G(t)>0$, $t\neq0,$ for some  $\mu\in(4,22^*)$, where $G(t)=\int^t_0 g(s)ds$,
\vskip 0.1 true cm
\noindent Yang et al. obtained a least energy sign-changing solution by means of the method of Nehari manifold, deformation arguments and
$L^\infty$-estimates. Here, we also would like to mention a more general quasilinear equation than equation \eqref{1.0.0}. Liu et al. \cite{LLW2} investigated the following quasilinear equation of the form
\begin{equation}\label{qlew}
\sum_{i, j=1}^{N} D_{j}\left(a_{i j}(x, u) D_{i} u\right)-\frac{1}{2} \sum_{i, j=1}^{N} D_{s} a_{i j}(x, u) D_{i} u D_{j} u+f(x, u)=0,
\end{equation}
where $D_{i}=\frac{\partial}{\partial x_{i}}$ and $D_{s} a_{i j}(x, s)=\frac{\partial}{\partial s} a_{i j}(x, s) .$  Equation \eqref{1.0.0} can be regarded as a special case of equation \eqref{qlew} for $a_{i j}(x, u)=\left(1+u^{2}\right) \delta_{i j}$.
By introducing a $p$-Laplacian perturbation approach, the authors obtained, via the method of invariant sets of descending flow, infinitely many sign-changing solutions of equation \eqref{qlew} in a bounded domain $\Omega\subset\R^N$. We point out that in \cite{LLW2} the following assumption is imposed: $\lim\limits _{s \rightarrow \infty} f(x, s) / s=+\infty$ uniformly with respect to $x \in \bar{\Omega}$ and for some $q>4$ and $c_0\in\R$,
$$
\frac{1}{q} s f(x, s)-F(x, s) \geq-c_{0}, \quad \forall x \in \bar{\Omega}, s \in \mathbb{R},
$$
where $F(x, s)=\int_{0}^{s} f(x, t) d t$. One can get that $\lim\limits _{s \rightarrow \infty} f(x, s) / s^3=+\infty$ uniformly with respect to $x \in \bar{\Omega}$.

Observe that the previous results on sign-changing solutions of (\ref{1.0.0}) were focused on the case where $g(u)=|u|^{p-2}u$ with $p\geq4$. A natural problem is whether
or not equation (\ref{1.0.0}) has sign-changing solutions if $p\in (2,4)$. In addition, for the case $p=4$, in \cite{LWW} Liu et al. considered the existence of a least energy sign-changing solution under the condition (V$'_1$), another interesting problem is whether or not a least energy sign-changing solution or higher
energy sign-changing solutions exist if (V$'_1$) is not satisfied?
In this paper, we will fill these gaps and give affirmative answers. Firstly we state the result about sign-changing solutions of (\ref{1.0.0}) with sub-cubic nonlinearity  and give the assumptions on $V$ and $g$ as follows.
  \vskip 0.2 true cm
 \noindent($V_1$) $V\in C(\mathbb{R}^{3})\cap L^\infty(\mathbb{R}^{3})$, $V(x)=V(|x|)$ and $\inf_{x\in\mathbb{R}^3}V(x)>0$.
\vskip 0.2 true cm \noindent($g_1$) $g\in C(\mathbb{R})$ and $g(t)=o(t)$\ as\ $t\rightarrow0$.
\vskip 0.2 true cm \noindent{($g_2$)} $\limsup_{|t|\rightarrow+\infty}\frac{|g(t)|}{|t|^{p-1}}<+\infty$ for some $p\in(2,12)$.
\vskip 0.2 true cm \noindent{($g_3$)}\ There exists $4\geq\mu>2$ such that $ g(t)t\geq \mu G(t)>0$ for $t\neq0$, where $G(t)=\int^t_0g(s)ds$.
\vskip 0.1 true cm
As a consequence of ($g_2$) and ($g_3$), one has $2<\mu\leq p<12$. In order to deal with sub-cubic nonlinearity, we also introduce the following condition on $V$.
\vskip 0.1 true cm
\noindent($V_2$) $V$ is weakly differentiable, $(\nabla V(x), x)\in L^\infty(\mathbb{R}^3)\cup L^{\frac32}(\mathbb{R}^3)$, and
$$\frac{\mu-2}{\mu}V(x)-(\nabla V(x), x)\geq0\ \ \ a.e. \ x\in\mathbb{R}^3,$$
where $\mu$ is given in ($g_3$) and $(\cdot,\cdot)$ is the usual inner product in $\mathbb{R}^3$.
\vskip 0.1 true cm
Now we state our first result.
\vskip 0.1 true cm
\noindent\textbf{Theorem 1.1.} {\it Let ($V_1$), ($V_2$) and ($g_1$)-($g_3$) hold. Then problem (\ref{1.0.0}) has a least energy sign-changing solution. If additionally $g$ is odd, then equation (\ref{1.0.0}) admits  infinitely many sign-changing solutions.}
\vskip 0.1 true cm
The main difficulties in the proof of Theorem 1.1 are three-fold.

First, a typical way in showing the existence and multiplicity of sign-changing solutions of (\ref{1.0.0}) is to use the method of sign-changing Nehari manifold, which depends on the monotonicity condition
\begin{equation}\label{1.0.2}\frac{g(t)}{t^3}\ \text{is increasing in}\ (0,+\infty)\ \text{and decreasing in}\ (-\infty,0).\end{equation}
Moreover, the Ambrosetti-Rabinowitz type condition: for some $\mu>4$,
\begin{equation}\label{1.0.3}\ g(t)t\geq \mu G(t)>0,\ \ t\neq0 \end{equation} is assumed to ensure the boundedness of Palais-Smale ((PS) for short) sequences.

In our paper, the nonlinearity is allowed to be the form of $g(t)=|t|^{p-2}t$ with $p\in(2,4]$, which does not satisfy (\ref{1.0.2}), let alone (\ref{1.0.3}), so the standard variational methods cannot be used directly.  Inspired by \cite{LWZ}, we will use a perturbation method and the method of invariant sets of descending flow. In \cite{LWZ}, Liu et al. considered the existence and multiplicity of sign-changing solutions for the Schr\"{o}dinger-Poisson system
\begin{equation*}\aligned
\left\{ \begin{array}{lll}
-\Delta u+V(x)u+\phi u=|u|^{p-2}u,\ & \text{in}\quad \mathbb{R}^3,\\
-\Delta\phi=u^2,    \ & \text{in}\quad
\mathbb{R}^3,
\end{array}\right.\endaligned
\end{equation*}
where $p\in(3,4]$, and  they overcame the difficulty brought by the term  $|u|^{p-2}u$ with $p\in(3,4]$ by adding a higher order nonlinear term and  the coercive condition ($V'_3$). Motivated by \cite{LWZ}, we shall use a perturbation method by adding a higher order term $\beta|u|^{r-2}u$ with $\beta>0$ and $r>4$. However, different from  \cite{LWZ}, we do not impose the coercive condition ($V'_3$), which plays an important role in showing the boundedness of (PS) sequence in \cite{LWZ}. So we shall add another perturbed term
$\lambda\Bigl(\int_{\mathbb{R}^3}u^2dx\Bigr)^\alpha u$ with $\lambda,\alpha>0$ to show the boundedness of (PS) sequences. We remark that the perturbation method is different from that in \cite{ZL}, where the authors  added a coercive potential term
and a $4$-Laplacian operator, and obtained infinitely many sign-changing solutions of (\ref{1.5}) with $p\in(4,22^*)$.

 Second, we shall use the method of invariant sets of descending flow to solve the perturbed problem, and go back to the original problem via the Pohozaev equality as in \cite{LWZ}. However, in this process, there are a lot of additional difficulties caused by the change $f$, which is used to transform a quasilinear problem into a semilinear one. Especially, the change $f$ leads that it is not easy to find a similar auxiliary operator $T_{\lambda,\beta}$ as in \cite{LWZ}, which  plays a crucial role in
constructing invariants sets of descending flow. In addition, some verifications have to be involved in the change $f$. For instance, when
showing the properties of the invariant sets of descending flow, we need to estimate $\left\langle u-T_{\lambda,\beta}(u), \frac{f(u)}{f'(u)}\right\rangle$ rather than $\left\langle u-T_{\lambda,\beta}(u), u\right\rangle$ as in \cite{LWZ}. When going back to the original problem, we have to combine a Pohozaev equality with the equality $\langle I'_{\lambda,\beta}(u),\frac{f(u)}{f'(u)}\rangle=0$ rather than $\langle I'_{\lambda,\beta}(u),u\rangle=0$, where $I_{\lambda,\beta}$ is the associated functional of the perturbed problem. Hence, some new estimations and tricks are required.

 Third, in the process in showing the multiplicity of sign-changing solutions, we will use the arguments of the existence part of Theorem 1.1 and some critical points theorem about multiple sign-changing solutions in \cite{LLW}. Firstly, for the perturbed problem with
two perturbation terms $\lambda\Bigl(\int_{\mathbb{R}^3}u^2dx\Bigr)^\alpha u$ and $\beta|u|^{r-2}u$, we
obtain infinitely many sign-changing solutions $u^j_{\lambda, \beta}$ with the energy $c^j_{\lambda, \beta}$, $j=1,2...$, and $c^j_{\lambda, \beta}\rightarrow+\infty$ as $j\rightarrow+\infty$. Then by taking $\lambda\rightarrow0^+$ and $\beta\rightarrow0^+$, sign-changing solutions $u^*_j$ with the energy $c^*_j$ of the original problem are obtained. To illustrate that there are infinitely many sign-changing solutions, we need to show $c^*_j\rightarrow+\infty$ as $j\rightarrow+\infty$. In \cite{LWZ}, there is just one perturbation term $\beta|u|^{r-2}u$, so it is easy to see that $c^*_j\geq c^j_\beta\rightarrow+\infty$, where $c^j_\beta,j=1,2,\cdots$ are the energy of sign-changing solutions of the perturbed problem in \cite{LWZ}.  In this paper, since the minimax values $c^j_{\lambda, \beta}$ have different monotonicity properties on the two
perturbation terms, it seems difficult to show $c^*_j\rightarrow+\infty$ by $c^j_{\lambda, \beta}\rightarrow+\infty$. As in \cite{LZ,LS}, by using
an auxiliary functional, we obtain $c^*_j\rightarrow+\infty$ as $j\rightarrow\infty$.

\begin{Remark} In this paper, we only consider the case $4\geq\mu>2$ in ($g_3$) since the case $\mu>4$ can be treated without any perturbation.
Actually, let ($V_1$), ($g_1$), ($g_2$) and ($g_3$) with $\mu>4$  hold, then equation (\ref{1.0.0}) admits a least energy sign-changing solution. If additionally $g$ is odd, then equation (\ref{1.0.0}) admits  infinitely many sign-changing solutions.
  \end{Remark}

Next we further study the case $\mu=4$ in ($g_3$), some typical problems can be considered without the condition (V$_2$). That is to say, we can deal with these problems directly without using the perturbation method. Now we  give the assumptions on $g$ and state the second result as follows.
\vskip 0.1 true cm
\noindent($g_4$) $G(t)>0$ for all $t\neq0$ and $\lim_{|t|\rightarrow+\infty}\frac{G(t)}{t^4}=+\infty$, where $G(t)=\int^t_0g(s)ds$.
\vskip 0.1 true cm
\noindent($g_5$) There exists $\gamma\in[1,+\infty)$ such that
$$\mathcal{G}(s)\leq \gamma \mathcal{G}(\tau),\ \text{for all}\ 0\leq s\leq \tau\ \text{or}\ \tau\leq s\leq 0,$$
where $\mathcal{G}(t)=\frac14 g(t)t-G(t)$ for any $t\in\mathbb{R}$.
\vskip 0.1 true cm

\begin{Remark} The condition ($g_5$) was firstly introduced by Jeanjean in \cite{Jeanjean}. ($g_5$) is weaker than (\ref{1.0.2}) since ($g_5$) holds with $\gamma=1$ if (\ref{1.0.2}) is satisfied.  By ($g_5$), it is easy to verify
\begin{equation}\label{1.4}\frac14 g(t)t-G(t)=\mathcal{G}(t)\geq0, \ \forall t\in\mathbb{R}.\end{equation} An example of $g$ satisfying ($g_1$), ($g_2$), ($g_4$) and ($g_5$) is $g(t)=t^3\log(1+|t|)$ for any $t\in\mathbb{R}$.
\end{Remark}

\noindent\textbf{Theorem 1.2.} {\it Let (V$_1$), ($g_1$), ($g_2$), ($g_4$) and ($g_5$) hold. Then problem (\ref{1.0.0}) has a least energy sign-changing solution. If additionally $g$ is odd, then equation (\ref{1.0.0}) admits  infinitely many sign-changing solutions.}

\begin{Remark}The proof of
Theorem 1.2 is based on the method of invariant sets of descending flow in \cite{LWZ} without any perturbation. However,
under the assumptions of Theorem 1.2, it seems difficult to verify
the associated functional satisfies the (PS) condition, which plays an important role in using the method of invariant sets of descending flow. In the present paper, we consider the Cerami ((Ce) for short) condition instead of the (PS) condition. So we have to establish a deformation lemma under the (Ce) condition and show that the associated functional satisfies the (Ce) condition.
 \end{Remark}
 \begin{Remark}
 We would like to point out that, if the functional satisfies the (Ce) condition under suitable assumptions on $V$ and $g$, then arguing as in Theorem 1.2, we can show the existence and multiplicity of sign-changing solutions for equation (\ref{1.0.0}). Moreover, we highlight that the method of invariant sets of descending flow under the (Ce) condition is also applicable for other related problems.
\end{Remark}

As a by-product, we finally state the result in the most typical case of $\mu=4$ in ($g_3$), i.e. $g(u)=u^3$, and consider
\begin{equation}\label{1.0.1}-\Delta u+V(x)u-u\Delta(u^2)=u^3,\ \ x\in\mathbb{R}^3.\end{equation}
\noindent\textbf{Theorem 1.3.} {\it Let (V$_1$) hold. Then problem (\ref{1.0.1}) has a least energy sign-changing solution. Moreover, equation (\ref{1.0.1}) admits  infinitely many sign-changing solutions.}
\begin{Remark}In \cite{LWW}, under the condition (V$'_1$), Liu et al. considered the existence of a least energy sign-changing solution of (\ref{1.0.1}) and there are few results about the multiplicity of sign-changing solutions. Hence, to  a certain extent, Theorem 1.3 completes the study made in \cite{LWW} on problem (\ref{1.0.1}).\end{Remark}

The paper is organized as
follows. In Section 2, we  give some preliminaries. In Sections 3, 4 and 5, we prove Theorems
1.1, 1.2 and 1.3 respectively.

\section{Preliminaries}
\renewcommand{\theequation}{2.\arabic{equation}}
In this paper we use the following
notations.  $\int_{\mathbb{R}^3}
f(x)dx$ is represented by $\int_{\mathbb{R}^3} f(x)$. For $2\leq p\leq\infty$, the
norm in $L^p(\mathbb{R}^3)$ is denoted by $|\cdot|_{p}$. For any $r>0$ and $x\in\mathbb{R}^3$,
$B_r(x)$ denotes the ball  centered at $x$ with the radius $r$.
The Hilbert space
$$E:= H^1_r(\mathbb{R}^3)=\bigl\{u\in H^1(\mathbb{R}^3):u(x)=u(|x|)\bigr\}$$ is the  space endowed with the following inner product and norm
$$(u,v)=\int_{\mathbb{R}^3} (\nabla u\nabla v+V(x)uv),\quad \|u\|^2=\int_{\mathbb{R}^3} (|\nabla u|^2+V(x)u^2).$$
Obviously, under (V$_1$), the norm $\|\cdot\|$ is an equivalent norm to the standard norm $\|\cdot\|_0=\bigl(|\nabla u|^2_2+|u|^2_2\bigr)^{\frac12}$ in $E$,  and so the embedding $E\hookrightarrow L^p(\mathbb{R}^3)$ is compact for any $p\in (2,6)$.

Due to the quasilinear term $-u\Delta(u^2)$, the functional of problem (\ref{1.0.0})
$$\mathcal{J}{(u)}=\frac{1}{2}\int_{\mathbb{R}^{3}} (1+2u^2)|\nabla u|^2+\frac 12\int_{\mathbb{R}^{3}}
V(x)u^2-\int_{\mathbb{R}^{3}} G(u),$$ is not well defined in $E$. As in
\cite{LWW}, we make use of a suitable change, namely, the change of
variables $v=f^{-1}(u)$, and then we can choose $E$
as the research space. The change $f$ is defined by $$\aligned
f'(t)&=\frac{1}{(1+2f^2(t))^{\frac12}} \ \text{on} \ [0,+\infty),\\
f(t)&=-f(-t) \ \ \ \ \ \ \text{on} \ (-\infty,0].\endaligned$$ Below we state some properties of $f$ given in \cite{CJ,DS1}.

\begin{lemma} \label{l2.1}
 \vskip 0.1 true cm \noindent(1) $f$ is uniquely defined, $C^\infty$
 and invertible;

\vskip 0.1 true cm \noindent (2) $|f'(t)|\leq1$ for all
$t\in\mathbb{R}$;

\vskip 0.1 true cm \noindent (3) $|f(t)|\leq|t|$ for all
$t\in\mathbb{R}$, and $|f(t)|\leq2^{\frac{1}{4}}|t|^{\frac12}$ for all $t\in\mathbb{R}$;

\vskip 0.1 true cm \noindent (4) $\frac{f(t)}{t}\rightarrow1$ as
$t\rightarrow0$;

\vskip 0.1 true cm \noindent (5)
$\frac{|f(t)|}{|t|^{\frac12}}\rightarrow2^{\frac{1}{4}}$ as
$|t|\rightarrow+\infty$;

\vskip 0.1 true cm \noindent (6) $\frac{f^2(t)}{2}\leq
tf'(t)f(t)\leq f^2(t)$ for all $t\in\mathbb{R}$;

\vskip 0.1 true cm \noindent (7)
there exists a positive constant $C$ such that
$|f(t)|\geq
C|t|$ if $|t|\leq1$, and
$|f(t)|\geq
C|t|^\frac12$ if  $|t|\geq1$;
\vskip 0.2 true cm \noindent (8) $|f(t)f'(t)|\leq\frac{1}{\sqrt{2}}$ for all $t\in\mathbb{R}$.
\end{lemma}

In view of the properties of the change $f$, from $\mathcal{{J}}$ we obtain
the functional
\begin{equation}\label{2.3}I(v)=\frac12\int_{\mathbb{R}^3}|\nabla v|^2+\frac12\int_{\mathbb{R}^3}
V(x)f^2(v)-\int_{\mathbb{R}^3} G(f(v)),\end{equation} which is well defined in
$E$ and of $C^1$ under our hypotheses. Moreover, the
critical points of $I$ are the weak solutions of the
problem\begin{equation}\label{2.0.0}
-\Delta v+V(x)f(v)f'(v)=g(f(v))f'(v),\ \ x\in\mathbb{R}^3.
\end{equation}
By the one-to-one correspondence $f$, we just need to study equation (\ref{2.0.0}). Denote
\begin{equation}\label{2.1}\tilde{g}(x,u)=g(f(u))f'(u)-V(x)f(u)f'(u)+V(x)u,\quad \tilde{G}(x,u)=\int^u_0\tilde{g}(x,s)ds.\end{equation}
One easily has the following properties of $\tilde{g}$.
\begin{lemma}\label{l3.1.0} For $q=\max\{3,\frac{p}{2}\}\in (2,6)$, there holds\\
\noindent(1) $\tilde{g}(x,t)=o(t)$ uniformly in $x$ as $t\rightarrow0$;
\vskip 0.2 true cm
\noindent(2) $\limsup_{|t|\rightarrow+\infty}\frac{|\tilde{g}(x,t)|}{|t|^{q-1}}<+\infty$ uniformly in $x$;\\
\noindent(3) for any $\epsilon>0$, there exists $C_\epsilon>0$ such that $|\tilde{g}(x,t)|\leq\epsilon|t|+C_\epsilon|t|^{q-1}$ for any $(x,t)\in\mathbb{R}^3\times\mathbb{R}$.
\end{lemma}

The following lemma plays an important role in showing the (PS) sequence is bounded.

\begin{lemma}\label{l2.3} Let (V$_1$) hold.
 Then there exists $C>0$ such that
 \begin{equation}\label{2.3.4}\|v\|\leq  C(|\nabla v|_2+|f(v)|^2_2+|f(v)|_2)\leq C\bigl(1+|\nabla v|_2+|f(v)|^2_2\bigr),\ \ \ \forall v\in E.\end{equation}
\end{lemma}
{\bf Proof}: For any $v\in E$, from Lemma \ref{l2.1} (7), the H\"{o}lder inequality and Young inequality it follows that
\begin{equation*}\aligned\int_{\{|v|>1\}}v^2&= \int_{\{|v|>1\}}|v|^{\frac65}|v|^{\frac45}\leq C\int_{\{|v|>1\}}|v|^{\frac65}|f(v)|^{\frac85}\\
&\leq C\Bigl(\int_{\{|v|>1\}}|v|^{6}\Bigr)^\frac15
\Bigl(\int_{\{|v|>1\}}|f(v)|^{2}\Bigr)^\frac45
\leq C\bigl(\frac35|v|^2_6+\frac25|f(v)|^4_2\bigr).\endaligned\end{equation*}
Using Lemma \ref{l2.1} (7) again we get
\begin{equation*}\int_{\{|v|\leq1\}}v^2\leq C\int_{\{|v|\leq1\}}f^2(v).\end{equation*}
Then
\begin{equation*}|v|_2\leq C(|v|_6+|f(v)|^2_2+|f(v)|_2)\leq C(1+|\nabla v|_2+|f(v)|^2_2).\end{equation*}
Therefore,
(\ref{2.3.4}) yields.\ \ \ $\Box$

In what follows, we recall the abstract critical point theorems developed by Liu et al. in \cite{LLW}, which will be used to show the existence and multiplicity of sign-changing solutions. Let $X$ be a Banach space, $J\in C^1(X,\mathbb{R})$ and $P, Q\subset X$ be open sets, $M=P\cap Q$, $\Sigma=\partial P\cap \partial Q$ and $W=P\cup Q$. For $c\in\mathbb{R}$, $K_c=\{x\in X: J(x)=c, \ J'(x)=0\}$ and $J^c=\{x\in X: J(x)\leq c\}$.

\begin{Definition}(\cite{LLW}) $\{P,Q\}$ is called an admissible family of invariant sets with respect to
$J$ at level $c$ provided that the following deformation property holds:
if $K_c\backslash W=\emptyset$, then there exists $\epsilon_0>0$ such that for $\epsilon\in (0,\epsilon_0)$, there exists $\eta\in C(X,X)$ satisfying

\noindent(1) $\eta(\overline{P})\subset\overline{P}$, $\eta(\overline{Q})\subset\overline{Q}$;\\
\noindent(2) $\eta|_{J^{c-2\epsilon}}=id$;\\
\noindent(3) $\eta(J^{c+\epsilon}\backslash W)\subset J^{c-\epsilon}$.
\end{Definition}

\begin{Theorem}\label{th1} (\cite{LLW}) Assume that $\{P, Q\}$ is an admissible family of invariant sets with respect
to $J$ at any level $c\geq c_*:=\inf_{u\in \Sigma}J(u)$ and there exists a map $\varphi_0:\triangle\rightarrow X$ satisfying

\noindent (1) $\varphi_0(\partial_1\triangle)\subset P$ and $\varphi_0(\partial_2\triangle)\subset Q$,\\
\noindent (2) $\varphi_0(\partial_0\triangle)\cap M=\emptyset$;\\
\noindent (3) $\sup_{u\in \varphi_0(\partial_0\triangle)}J(u)<c_*$,\\
where $\triangle=\{(t_1,t_2)\in\mathbb{R}^2:t_1,t_2\geq0, t_1+t_2\leq1\}$, $\partial_1\triangle=\{0\}\times[0,1]$, $\partial_2\triangle=[0,1]\times\{0\}$ and $\partial_0\triangle=\{(t_1,t_2)\in\mathbb{R}^2:t_1,t_2\geq0, t_1+t_2=1\}$. Define
$$c=\inf_{\varphi\in \Gamma}\sup_{u\in \varphi(\triangle)\backslash W}J(u),$$
where $\Gamma:=\{\varphi\in C(\triangle, E):\varphi(\partial_1\triangle)\subset P, \varphi(\partial_2\triangle)\subset Q, \varphi|_{\partial_0\triangle}={\varphi_0}|_{\partial_0\triangle}\}$.
Then $c\geq c_*$ and $K_c\backslash W\neq\emptyset$.
\end{Theorem}

 If additionally assume $G: X\rightarrow X$ is an isometric involution, i.e. $G^2=id$ and $d(Gx,Gy)=d(x,y)$ for $x,y\in X$. We assume $J$ is $G-$invariant on $X$ in the sense that $J(Gx)=J(x)$ for any $x\in X$. We also assume $Q=GP$. A subset $F\subset X$ is said to be symmetric if $Gx\in F$ for any $x\in F$. The genus of a closed symmetric subset $F$ of $X\backslash\{0\}$ is denoted by $\gamma(F)$.

\begin{Definition}(\cite{LLW}) $P$ is called a $G-$admissible invariant set with respect to $J$ at
level $c$, if the following deformation property holds: there exist $\epsilon_0>0$ and a symmetric
open neighborhood $N$ of $K_c\backslash W$ with $\gamma(\overline{N})<+\infty$, such that for $\epsilon\in (0,\epsilon_0)$ there exists $\eta\in C(X,X)$ satisfying\\
\noindent(1) $\eta(\overline{P})\subset\overline{P}$, $\eta(\overline{Q})\subset\overline{Q}$;\\
\noindent(2) $\eta\circ G=G\circ\eta$;\\
\noindent(3)$\eta|_{J^{c-2\epsilon}}=id$;\\
\noindent(4)$\eta(J^{c+\epsilon}\backslash(N\cup W))\subset J^{c-\epsilon}$.\end{Definition}

\begin{Theorem}\label{th2} (\cite{LLW}) Assume that $P$ is a $G-$admissible invariant set with respect to $J$ at
level $c\geq c_*:=\inf_{u\in \Sigma}J(u)$ and for any $n\in\mathbb{N}$, there exists a continuous map $\varphi_n:B_n=\{x\in \mathbb{R}^n:|x|\leq1\}\rightarrow X$ satisfying\\
\noindent(1) $\varphi_n(0)\in M:=P\cap Q$, $\varphi_n(-t)=G\varphi_n(t)$ for $t\in B_n$;\\
\noindent(2) $\varphi_n(\partial B_n)\cap M=\emptyset$;\\
\noindent(3) $\sup_{u\in {Fix_G\cup\varphi_n(\partial B_n)}}J(u)<c_*$, where $Fix_G:=\{u\in X:Gu=u\}$. \\
For $j\in\mathbb{N}$, define
$$c_j=\inf_{B\in \Gamma_j}\sup_{u\in B\backslash W}J(u),$$
where
$$\aligned\Gamma_j:=\{B|B=&\varphi(B_n\backslash Y)\ \text{for some}\ \varphi\in G_n, n\geq j,\\&\ \text{and open}\ Y\subset B_n\ \text{such that}\ -Y=Y\ \text{and}\ \gamma(\overline{Y})\leq n-j\}.\endaligned$$
and
$$G_n:=\{\varphi|\varphi\in C(B_n, X), \varphi(-t)=G\varphi(t) \ \text{for}\ t\in B_n, \varphi(0)\in M\ \text{and}\ \varphi|_{\partial B_n}={\varphi_n}|_{\partial B_n}\}.$$
Then for $j\geq2$, $c_j\geq c_*$, $K_{c_j}\backslash W\neq\emptyset$ and $c_j\rightarrow+\infty$ as $j\rightarrow+\infty$.
 \end{Theorem}

\section{ Proof of Theorem 1.1}
\renewcommand{\theequation}{3.\arabic{equation}}
In this section, we will show Theorem 1.1 and assume ($V_1$), ($V_2$) and ($g_1$)-($g_3$) are satisfied. Since $\mu\leq4$ in ($g_3$), it is not easy to show the (PS) sequence is bounded. A perturbed problem is introduced to overcome this difficulty.
\subsection{A perturbed problem}
 Set $\alpha\in (0,\frac{\mu-2}{3\mu+2})$  and fix $\lambda, \beta\in (0,1]$ and $r\in(\max\{4,p\},12)$, we consider the perturbed problem of (\ref{2.0.0}) that
\begin{equation}\label{3.0.0}-\Delta u+V(x)f(u)f'(u)+\lambda|f(u)|^{2\alpha}_2 f(u)f'(u)=g(f(u))f'(u)+\beta|f(u)|^{r-2}f(u)f'(u).\end{equation}
The associated functional is
$$I_{\lambda,\beta}(u)=I(u)+
\frac{\lambda}{2(1+\alpha)}|f(u)|^{2(1+\alpha)}_2
-\frac{\beta}{r}|f(u)|^r_r,$$
where $I$ is given in (\ref{2.3}).
In order to go back to the original problem (\ref{2.0.0}),
we will make use of the Pohozaev type identity of (\ref{3.0.0}), whose proof is standard and can be referred to \cite[Lemma 2.6]{XC} for example.
\begin{lemma} \label{l4.1} Let $u$ be a critical point of $I_{\lambda,\beta}$ for $(\lambda,\beta)\in(0,1]\times(0,1]$, then
$$\aligned \frac12|\nabla u|^2_2&+\frac12\int_{\mathbb{R}^3}[3V(x)+(\nabla V(x), x)]f^2(u)\\&+\frac{3\lambda}{2(1+\alpha)}|f(u)|^{2(1+\alpha)}_2
-3\int_{\mathbb{R}^3}G(f(u))-\frac{3\beta}{r}|f(u)|^r_r=0.\endaligned$$\end{lemma}

We now introduce an auxiliary operator $T_{\lambda,\beta}$, which will be used to construct the descending flow for
the functional $I_{\lambda,\beta}$. As  an application of Lax-Milgram theorem, for any $u\in E$, there is a unique solution $v=T_{\lambda,\beta}(u)\in E$
of the equation
\begin{equation}\label{3.1.0} -\Delta v+V(x)v+\lambda|f(u)|^{2\alpha}_2\frac{f(u)f'(u)}{u}v=\bar{g}(x,u),\end{equation}
where $$\bar{{g}}(x,u)=g(f(u))f'(u)+\beta|f(u)|^{r-2}f(u)f'(u)-
V(x)f(u)f'(u)+V(x)u,$$
  and $\bar{G}(x,u):=\int^u_0\bar{{g}}(x,s)ds$. Clearly, the three statements are equivalent: $u$ is a solution of (\ref{3.0.0}), $u$ is a critical point of $I_{\lambda,\beta}$,
and $u$ is a fixed point of $T_{\lambda,\beta}$. Moreover, $\bar{g}$ has the following properties.
 \begin{lemma}\label{l4.1.0} (1) $\bar{g}(x,t)=o(t)$ uniformly in $x$ as $t\rightarrow0$;
\vskip 0.1 true cm
 \noindent (2) for any $\epsilon>0$, there exists $C_\epsilon>0$ such that
\begin{equation}\label{4.0.2} |\bar{{g}}(x,t)|\leq\epsilon|t|+C_\epsilon|t|^{\frac r2-1},\quad \forall (x,t)\in\mathbb{R}^3\times\mathbb{R}. \end{equation}
\end{lemma}

\begin{lemma}\label{l4.2}$T_{\lambda,\beta}$ is continuous and compact.\end{lemma}
{\bf Proof}: We firstly show that $T_{\lambda,\beta}$ is continuous.  Assume that $u_n\rightarrow u$ in $E$. Up to a subsequence, suppose that $u_n\rightarrow u$ in $L^s(\mathbb{R}^3)$ with $s\in (2,6)$. Set $v_n=T_{\lambda,\beta}(u_n)$ and $v=T_{\lambda,\beta}(u)$, we have
\begin{equation}\label{4.1.2}-\Delta v_n+V(x)v_n+\lambda|f(u_n)|^{2\alpha}_2\frac{f(u_n)f'(u_n)}{u_n}v_n
=\bar{{g}}(x,u_n),\end{equation}
and \begin{equation}\label{4.1.3}-\Delta v+V(x)v+\lambda|f(u)|^{2\alpha}_2\frac{f(u)f'(u)}{u}v
=\bar{g}(x,u).\end{equation}
Testing with $v_n$ in (\ref{4.1.2}), by (\ref{4.0.2}) and $r\in(4,12)$ we get
\begin{equation}\label{4.2.2}\aligned\|v_n\|^2+\lambda|f(u_n)|^{2\alpha}_2\int_{\mathbb{R}^3}\frac{f(u_n)f'(u_n)}{u_n}v^2_n
&=\int_{\mathbb{R}^3}\bar{g}(x,u_n)v_n\\
&\leq\epsilon\|u_n\|\|v_n\|+C_\epsilon\|u_n\|^{{\frac{r}{2}-1}}\|v_n\|.\endaligned\end{equation}
Then $\{v_n\}$ is bounded in $E$. After passing to a subsequence, suppose $v_n\rightharpoonup v^*$ in $E$, $v_n\rightarrow v^*$ in $L^2_{loc}(\mathbb{R}^3)$ and $v_n\rightarrow v^*$ in $L^s(\mathbb{R}^3)$ with $s\in(2,6)$. Using (\ref{4.1.2}) it is easy to see that $v^*$ is a solution of (\ref{4.1.3}) and then $v^*=v$ using the uniqueness. Moreover, testing with $v_n-v$ in (\ref{4.1.2}) and (\ref{4.1.3}) we obtain
\begin{equation}\label{3.6.0}\|v_n-v\|^2+\lambda\mathbb{A}_1-\lambda\mathbb{A}_2=
\int_{\mathbb{R}^3}(\bar{g}(x,u_n)
-\bar{g}(x,u))(v_n-v),\end{equation}
where $$\mathbb{A}_1:=|f(u_n)|^{2\alpha}_2\int_{\mathbb{R}^3}\frac{f(u_n)f'(u_n)}{u_n}v_n(v_n-v),
\ \mathbb{A}_2:=|f(u)|^{2\alpha}_2\int_{\mathbb{R}^3}\frac{f(u)f'(u)}{u}v(v_n-v).$$
Since $v_n\rightarrow v$ in $L^2_{loc}(\mathbb{R}^3)$, $|f(u_n)|_2$ is bounded and the fact that $\Bigl|\frac{f(t)f'(t)}{t}\Bigr|\leq1$ for any $t\neq0$, one easily has  $$\aligned
&|f(u_n)|^{2\alpha}_2\Bigl|\int_{\mathbb{R}^3}\Bigl[\frac{f(u_n)f'(u_n)}{u_n}-\frac{f(u)f'(u)}{u}\Bigr]v(v_n-v)\Bigr|
\leq2|f(u_n)|^{2\alpha}_2\int_{\mathbb{R}^3}|v(v_n-v)|=o_n(1),\\
&\Bigl(|f(u)|^{2\alpha}_2-|f(u_n)|^{2\alpha}_2\Bigr)\int_{\mathbb{R}^3}\frac{f(u)f'(u)}{u}v(v_n-v)=o_n(1).\endaligned$$
 Then
\begin{equation}\label{4.2.0}\aligned\mathbb{A}_1-\mathbb{A}_2=&
|f(u_n)|^{2\alpha}_2\int_{\mathbb{R}^3}\frac{f(u_n)f'(u_n)}{u_n}(v_n-v)^2\\&+
|f(u_n)|^{2\alpha}_2\int_{\mathbb{R}^3}\Bigl[\frac{f(u_n)f'(u_n)}{u_n}-\frac{f(u)f'(u)}{u}\Bigr]v(v_n-v)
\\&+\Bigl(|f(u_n)|^{2\alpha}_2-|f(u)|^{2\alpha}_2\Bigr)\int_{\mathbb{R}^3}\frac{f(u)f'(u)}{u}v(v_n-v)\\
=&|f(u_n)|^{2\alpha}_2\int_{\mathbb{R}^3}\frac{f(u_n)f'(u_n)}{u_n}(v_n-v)^2+o_n(1).\endaligned\end{equation}
Next we show \begin{equation}\label{4.2.1}\mathbb{B}_n:=\int_{\mathbb{R}^3}(\bar{g}(x,u_n)
-\bar{g}(x,u))(v_n-v)\rightarrow0.\end{equation} In fact, let $\phi\in C^\infty_0(\mathbb{R},[0,1])$ be such that $\phi(t)=1$ for $|t|\leq1$ and $\phi(t)=0$ for $|t|\geq2$. Setting
$$\bar{{g}}_1(x,t)=\phi(t)\bar{{g}}_1(x,t),\quad \quad \bar{{g}}_2(x,t)=\bar{{g}}(x,t)-\bar{{g}}_1(x,t).$$
By (\ref{4.0.2}), there exists $C>0$ such that
$$|\bar{{g}}_1(x,t)|\leq C|t|,\ |\bar{{g}}_2(x,t)|\leq C|t|^{\frac{r}{2}-1},\ \ \text{for all}\ (x,t)\in\mathbb{R}^3\times\mathbb{R}.$$
Therefore
$$\aligned
\mathbb{B}_n&=\int_{\mathbb{R}^3}(\bar{{g}}_1(x,u_n)-\bar{{g}}_1(x,u))(v_n-v)+
\int_{\mathbb{R}^3}(\bar{{g}}_2(x,u_n)-\bar{{g}}_2(x,u))(v_n-v)
\\&\leq \epsilon\int_{\mathbb{R}^3}(|u_n|+|u|)|v_n-v|+
|\bar{{g}}_2(x,u_n)-\bar{{g}}_2(x,u)|_{\frac{r}{r-2}}
|v_n-v|_{\frac{r}{2}}\leq C\epsilon+o_n(1).\endaligned$$
By the arbitrariness of $\epsilon$, we know $\mathbb{B}_n\rightarrow0$ as $n\rightarrow\infty$. Combining with (\ref{3.6.0}), (\ref{4.2.0}) and the fact that $\frac{f(t)f'(t)}{t}\geq0$ for any $t\neq0$, we have $v_n\rightarrow v$ in $E$.

Below we show that $T_{\lambda,\beta}$ is compact. Assume that $u_n\rightharpoonup u$ in $E$ and $u_n\rightarrow u$ in $L^s(\mathbb{R}^3)$ with $s\in (2,6)$. Set $\lim_{n\rightarrow\infty}|f(u_n)|^{2\alpha}_2=a$ and $v_n=T_{\lambda,\beta}(u_n)$. As the above argument we infer $\{v_n\}$ is bounded in $E$. Up to a subsequence, suppose that $v_n\rightharpoonup v_0$ in $E$, $v_n\rightarrow v_0$ in $L^s(\mathbb{R}^3)$ with $s\in(2,6)$ and $v_n(x)\rightarrow v_0(x)$ a.e. in $\mathbb{R}^3$. For any $w\in C^\infty_0(\mathbb{R}^3)$,
by the fact that $|\frac{f(t)f'(t)}{t}|\leq1$ for any $t\neq0$, we know
\begin{equation*}\int_{\mathbb{R}^3}\frac{f(u_n)f'(u_n)}{u_n}(v_n-v_0)w
=o_n(1),\end{equation*}
and using the mean value theorem we infer
\begin{equation*}\aligned
&\Bigl|\int_{\mathbb{R}^3}\Bigl[\frac{f(u_n)f'(u_n)}{u_n}
-\frac{f(u)f'(u)}{u}\Bigr]v_0w\Bigr|=
\Bigl|\int_{\mathbb{R}^3}\frac{f(\xi_n)f'(\xi_n)}{\xi_n}(u_n-u)v_0w\Bigr|
\\ \leq&\int_{\mathbb{R}^3}\bigl|(u_n-u)v_0w\bigr|\leq |u_n-u|_3
|v_0|_3|w|_3=o_n(1),\endaligned\end{equation*}
where $\xi_n=\theta u_n+(1-\theta)u$ for some $\theta\in(0,1)$.
Then there holds
 $$\aligned &|f(u_n)|^{2\alpha}_2\int_{\mathbb{R}^3}\frac{f(u_n)f'(u_n)}{u_n}v_nw-
a\int_{\mathbb{R}^3}\frac{f(u)f'(u)}{u}v_0w\\=&a\int_{\mathbb{R}^3}\frac{f(u_n)f'(u_n)}{u_n}v_nw-
a\int_{\mathbb{R}^3}\frac{f(u)f'(u)}{u}v_0w+o_n(1)\\=&
a\int_{\mathbb{R}^3}\frac{f(u_n)f'(u_n)}{u_n}(v_n-v_0)w
+a\int_{\mathbb{R}^3}\Bigl[\frac{f(u_n)f'(u_n)}{u_n}
-\frac{f(u)f'(u)}{u}\Bigr]v_0w=o_n(1).
\endaligned$$
In view of (\ref{4.1.2}) we know $v_0$ is a solution of the equation \begin{equation*}-\Delta v+V(x)v+\lambda a\frac{f(u)f'(u)}{u}v
=\bar{g}(x,u).\end{equation*}
Then
\begin{equation}\label{4.2.3}\|v_0\|^2+\lambda a\int_{\mathbb{R}^3}\frac{f(u)f'(u)}{u}v^2_0=\int_{\mathbb{R}^3}\bar{g}(x,u)v_0.\end{equation}
In the same way as (\ref{4.2.1}), we have
$$\int_{\mathbb{R}^3}(\bar{g}(x,u_n)-\bar{g}(x,u))v_n=o_n(1).$$
and so
\begin{equation}\label{4.2.7}\int_{\mathbb{R}^3}\bigl(\bar{g}(x,u_n)v_n-\bar{g}(x,u)v_0\bigr)=o_n(1)+
\int_{\mathbb{R}^3}\bar{g}(x,u)(v_n-v_0)=o_n(1).\end{equation}
From Fatou lemma, the equality in (\ref{4.2.2}) and (\ref{4.2.7}) it follows that
$$\aligned\|v_0\|^2+\lambda a\int_{\mathbb{R}^3}\frac{f(u)f'(u)}{u}v^2_0&
\leq\liminf_{n\rightarrow\infty}\Bigl[\|v_n\|^2+\lambda|f(u_n)|^{2\alpha}_2
\int_{\mathbb{R}^3}\frac{f(u_n)f'(u_n)}{u_n}v^2_n\Bigr]\\
&=\liminf_{n\rightarrow\infty}\int_{\mathbb{R}^3}\bar{g}(x,u_n)v_n=\int_{\mathbb{R}^3}\bar{g}(x,u)v_0.\endaligned$$
Using (\ref{4.2.3}) we know $v_n\rightarrow v_0$ in $E$. This ends the proof.
 \ \ \ \ \ \ $\Box$

\begin{lemma}\label{l4.3}(1) $\langle I'_{\lambda,\beta}(u), u-T_{\lambda,\beta}(u)\rangle\geq\|u-T_{\lambda,\beta}(u)\|^2$ for all $u\in E$;\\
(2) $\|I'_{\lambda,\beta}(u)\|\leq\|u-T_{\lambda,\beta}(u)\|(1+C_1\|u\|^{2\alpha})$ for all $u\in E$, where $C_1>0$ is a positive constant independent on $\lambda$ and $\beta$.  \end{lemma}
{\bf Proof}: (1) Since $T_{\lambda,\beta}(u)$ is the solution of (\ref{3.1.0}) and the fact that $\frac{f(t)f'(t)}{t}\geq0$ for any $t\neq0$, we obtain $$\aligned\langle I'_{\lambda,\beta}(u), u-T_{\lambda,\beta}(u)\rangle&
=\|u-T_{\lambda,\beta}(u)\|^2+\lambda|f(u)|^{2\alpha}_2
\int_{\mathbb{R}^3}\frac{f(u)f'(u)}{u}[u-T_{\lambda,\beta}(u)]^2
\\&\geq\|u-T_{\lambda,\beta}(u)\|^2.\endaligned$$

\noindent(2) For any $\varphi\in E$, by Lemma \ref{l2.1} (2) and (3) we infer
$$\aligned\langle I'_{\lambda,\beta}(u), \varphi\rangle&=\langle u-T_{\lambda,\beta}(u), \varphi\rangle+\lambda|f(u)|^{2\alpha}_2
\int_{\mathbb{R}^3}\frac{f(u)f'(u)}{u}(u-T_{\lambda,\beta}(u))\varphi\\
&\leq\|u-T_{\lambda,\beta}(u)\|\|\varphi\|+\lambda|u|^{2\alpha}_2\|u-T_{\lambda,\beta}(u)\|\|\varphi\|.\endaligned$$
Then $\|I'_{\lambda,\beta}(u)\|\leq\|u-T_{\lambda,\beta}(u)\|
(1+C_1\|u\|^{2\alpha})$.\ \ \ \ \ \ $\Box$

\begin{lemma}\label{l4.4.0} Fix $(\lambda,\beta)\in (0,1]\times(0,1]$ and assume $a<b$ and $\tau>0$.  If $u\in E$, $I_{\lambda,\beta}(u)\in [a,b]$ and $\|I'_{\lambda,\beta}(u)\|\geq\tau$, then there exists $\delta>0$ (which depends on $\lambda$ and $\beta$) such that $\|u-T_{\lambda,\beta}(u)\|\geq\delta$.\end{lemma}
{\bf Proof}: By Lemma \ref{l2.1} (6), for any $u\in E$ we have
\begin{equation*}\nabla \bigl({\frac{f(u)}{f'(u)}}\bigr)=\frac{1+4f^2(u)}{1+2f^2(u)}\nabla u,\ \Bigl|\frac{f(u)}{f'(u)}\Bigr|=\Bigl|\frac{f(u)u}{f'(u)u}\Bigr|\leq
\Bigl|\frac{f(u)u}{\frac{f(u)}2}\Bigr|
\leq  2|u|.\end{equation*}
Then \begin{equation}\label{3.0.1}\bigl\|\frac{f(u)}{f'(u)}\bigr\|\leq 2\|u\|.\end{equation}
Taking $\gamma\in (4,r)$,  we obtain
 $$\aligned &I_{\lambda,\beta}(u)-\frac{1}{\gamma}\langle u-T_{\lambda,\beta}(u),\frac{f(u)}{f'(u)}\rangle\\
 =&I_{\lambda,\beta}(u)-\frac{1}{\gamma}\langle I'_{\lambda,\beta}(u),\frac{f(u)}{f'(u)}\rangle+\frac{\lambda}{\gamma}|f(u)|^{2\alpha}_2\int_{\mathbb{R}^3}\frac{f^2(u)}{u}(u-T_{{\lambda,\beta}}(u))
\\=& \int_{\mathbb{R}^3}\Bigl(\frac12-\frac{1+4f^2(u)}{{\gamma}(1+2f^2(u))}\Bigr)|\nabla u|^2_2+\Bigl(\frac12-\frac1{\gamma}\Bigr)\int_{\mathbb{R}^3}V(x)f^2(u)\\&+\int_{\mathbb{R}^3}\Bigl[\frac1\gamma g(f(u))f(u)-G(f(u))\Bigr]
+\frac{(r-\gamma)\beta}{r\gamma}|f(u)|^r_r\\&+\Bigl(\frac{1}{2(1+\alpha)}-\frac{1}{\gamma}\Bigr)\lambda|f(u)|^{2(1+\alpha)}_2
+\frac{\lambda}{\gamma}|f(u)|^{2\alpha}_2 \int_{\mathbb{R}^3}\frac{f^2(u)}{u}(u-T_{{\lambda,\beta}}(u))
.\endaligned$$
Using ($g_1$) and ($g_2$), for any $\epsilon>0$, there exists a constant
$C_{\epsilon}>0$ such that
\begin{equation}\label{2.2}|g(t)|\leq \epsilon|t|+C_\epsilon|t|^{p-1},\ \forall t\in\mathbb{R}.\end{equation}
Then $$\aligned &I_{\lambda,\beta}(u)-\frac{1}{\gamma}\langle u-T_{\lambda,\beta}(u),\frac{f(u)}{f'(u)}\rangle\\ \geq&
\Bigl(\frac12-\frac2{\gamma}-\epsilon\Bigr)\bigl[|\nabla u|^2_2+\int_{\mathbb{R}^3}V(x)f^2(u)\bigr]-C_\epsilon|f(u)|^p_p+
\frac{(r-\gamma)\beta}{r\gamma}|f(u)|^r_r\\&+\Bigl(\frac{1}{2(1+\alpha)}
-\frac{1}{\gamma}\Bigr)\lambda|f(u)|^{2(1+\alpha)}_2-
\frac{\lambda}{\gamma}|f(u)|^{2\alpha}_2 |u|_2|u-T_{{\lambda,\beta}}(u)|_2.\endaligned$$
Choosing $\epsilon$ small enough, by (\ref{2.3.4}) and (\ref{3.0.1}) we deduce
\begin{equation}\label{4.7.0}\aligned &|\nabla u|^2_2+\int_{\mathbb{R}^3}V(x)f^2(u)
+\lambda|f(u)|^{2(1+\alpha)}_2+\beta|f(u)|^r_r-C_0|f(u)|^p_p
\\ \leq& C\Bigl[|I_{\lambda,\beta}(u)|+2\|u-T_{\lambda,\beta}(u)\|(|\nabla u|_2+1+|f(u)|^2_2)\\&
+|f(u)|^{2\alpha}_2(|\nabla u|_2+1+|f(u)|^2_2)\|u-T_{\lambda,\beta}(u)\|
\Bigr].\endaligned\end{equation}

 Now we assume on the contrary that there exists $\{u_n\}\subset E$ with $I_{\lambda,\beta}(u_n)\in[a,b]$ and $\|I'_{\lambda,\beta}(u_n)\|\geq \alpha$ such that
$$\|u_n-T_{\lambda,\beta}(u_n)\|\rightarrow0, \quad\text{as}\ n\rightarrow\infty.$$
By (\ref{4.7.0}) we get
\begin{equation}\label{4.7.1}\aligned &|\nabla u_n|^2_2+\int_{\mathbb{R}^3}V(x)f^2(u_n)
+{\lambda}|f(u_n)|^{2(1+\alpha)}_2+\beta|f(u_n)|^r_r-C_0|f(u_n)|^p_p
\\
\leq&C\Bigl[C_1+2\|u_n-T_{\lambda,\beta}(u_n)\|(|\nabla u_n|_2+1+|f(u_n)|^2_2)\\&
+(|f(u_n)|^{4\alpha}_2+|\nabla u_n|^2_2+|f(u_n)|^{2\alpha}_2+|f(u_n)|^{2+2\alpha}_2)\|u_n-T_{\lambda,\beta}(u_n)\|
\Bigr].\endaligned\end{equation}
In the same way as (3.7) in
 \cite{LZ}, for any $c,d>0$ and $2<p<r<12$, there holds
 \begin{equation}\label{4.1.5}\inf_{u\in H^1(\mathbb{R}^3)}(|f(u)|^r_r+c|f(u)|^{2\alpha+2}_2-d|f(u)|^p_p)>-\infty.\end{equation}
Therefore
$$\frac{\lambda}{2}|f(u_n)|^{2(1+\alpha)}_2+\beta|f(u_n)|^r_r
-C_{0}|f(u_n)|^p_p\geq -C_2,$$
where $C_2>0$ is independent on $n$.
From (\ref{4.7.1}) it follows that
$$\aligned &|\nabla u_n|^2_2+\int_{\mathbb{R}^3}V(x)f^2(u_n)+
\frac{\lambda}{2}|f(u_n)|^{2(1+\alpha)}_2\\
\leq &
C_2+C\Bigl[C_1+2\|u_n-T_{\lambda,\beta}(u_n)\|(|\nabla u_n|_2+1+|f(u_n)|^2_2)\\&
+\bigl(|f(u_n)|^{4\alpha}_2+|\nabla u_n|^2_2+|f(u_n)|^{2\alpha}_2+|f(u_n)|^{2+2\alpha}_2\bigr)\|u_n-T_{\lambda,\beta}(u_n)\|
\Bigr].
\endaligned$$
Since $\alpha<1$ we deduce that
 $\{|\nabla u_n|_2+|f(u_n)|_2\}$ {is bounded} and so $\{u_n\}$ is bounded in $E$ using (\ref{2.3.4}). According to Lemma \ref{l4.3} (2) we know $\|I'_{\lambda,\beta}(u_n)\|\rightarrow0$ as $n\rightarrow\infty$. This is impossible since $\|I'_{\lambda,\beta}(u_n)\|\geq\alpha$. The proof is complete.\ \ \ \ \ $\Box$

\subsection{Invariant subsets of descending flows}
To look for sign-changing solutions, we define the positive and negative cones by
$$P^+:=\{u\in E: u\geq0\} \ \text{and}\ P^-:=\{u\in E: u\leq0\}$$
respectively. For $\epsilon>0$, set
$$P^+_\epsilon:=\{u\in E: dist(u,P^+)<\epsilon\} \ \text{and}\ P^-_\epsilon:=\{u\in E: dist(u,P^+)<\epsilon\},$$
where $dist(u, P^{\pm})=\inf_{u\in P^{\pm}}\|u-v\|$. Clearly, $P^{-}_\epsilon=-P^{+}_\epsilon$.  Let $W=P^{+}_\epsilon\cup P^{-}_\epsilon$. Then, $W$
is an open and symmetric subset of $E$ and $E\backslash{W}$ contains only sign-changing functions. In the following lemma, we will show that, for $\epsilon$ small enough, all sign-changing solutions to (\ref{3.0.0})
are contained in $E\backslash{W}$.

\begin{lemma}\label{l4.4} There is $\epsilon_0>0$ independent on $\lambda$ and $\beta$ such that for $\epsilon\in(0,\epsilon_0)$,

\noindent(1) $T_{\lambda,\beta}(\partial P^{-}_\epsilon)\subset P^{-}_\epsilon$ and every nontrivial solution $u\in P^{-}_\epsilon$ is negative.

\noindent(2) $T_{\lambda,\beta}(\partial P^{+}_\epsilon)\subset P^{+}_\epsilon$ and every nontrivial solution $u\in P^{+}_\epsilon$ is positive.
\end{lemma}
{\bf Proof}:
We only prove (1) since the argument of (2) is similar. For $u\in E$, define $v:=T_{\lambda,\beta}(u)$. Note that
$$\frac{f(t)f'(t)}{t}\geq0,\ \text{for all}\ t\neq0.$$
Then using (\ref{4.0.2}), for any
$\delta>0$, there exists $C_\delta>0$ such that
$$\aligned dist(v, P^-)\|v^+\|&\leq\|v^+\|^2=\int_{\mathbb{R}^3}\bar{g}(x,u)v^+-\lambda|f(u)|^{2\alpha}_2\int_{\mathbb{R}^3}\frac{f(u)f'(u)}{u}vv^+
\\&\leq\int_{\mathbb{R}^3}\bar{g}(x,u^+)v^+\leq\int_{\mathbb{R}^3}(\delta|u^+|+C_\delta|u^+|^{\frac{r}{2}-1})v^+\\
&\leq C(\delta dist(u, P^-)+C_\delta dist(u,P^-)^{\frac{r}{2}-1})\|v^+\|.\endaligned$$
Thus
$$ dist(v, P^-)\leq C(\delta dist(u, P^-)+C_\delta dist(u,P^-)^{\frac r2-1}).$$
Choosing $\delta$ small enough, there exists $\epsilon_0>0$ such that for $\epsilon\in (0,\epsilon_0)$,
$$dist (T_{\lambda,\beta}(u), P^-)\leq\frac12dist(u, P^-),\quad \forall u\in P^-_\epsilon.$$
So $T_{\lambda,\beta}(\partial P^{-}_\epsilon)\subset P^{-}_\epsilon$. If there exists $u\in P^-_\epsilon$ such that $T_{\lambda,\beta}(u)=u$, then $u\in P^-$. If $u\not\equiv0$, then the maximum principle implies that $u<0$ in $\mathbb{R}^3$.\ \ \ \ \ $\Box$

Denote the set of fixed points of $T_{\lambda,\beta}$ by $K$. Since the operator $T_{\lambda,\beta}$ is not locally Lipschitz continuous, we need to construct a locally Lipschitz continuous operator $B_{\lambda,\beta}$ on $E_0:=E\backslash K$ which inherits its properties. Arguing as the
proof of \cite[Lemma 2.1]{BLW}, we have the following result.

\begin{lemma}\label{l4.5}There exists a locally Lipschitz continuous operator $B_{\lambda,\beta}: E_0\rightarrow E$ such
that \\
(1)\ $B_{\lambda,\beta}(\partial P^+_\epsilon)\subset P^+_\epsilon$ and $B_{\lambda,\beta}(\partial P^-_\epsilon)\subset P^-_\epsilon$ for $\epsilon\in (0,\epsilon_0)$;

\noindent(2)\ $\frac12\|u-B_{\lambda,\beta}(u)\|\leq\|u-T_{\lambda,\beta}(u)\|\leq 2\|u-B_{\lambda,\beta}(u)\|$ for all $u\in E_0$;

\noindent(3)\ $\langle I'_{\lambda,\beta}(u), u-B_{\lambda,\beta}(u)\rangle\geq\frac12\|u-T_{\lambda,\beta}(u)\|^2$ for all $u\in E_0$;

\noindent(4)\ if $g$ is odd, then $B_{\lambda,\beta}$ is odd.
\end{lemma}

   In the following, we show that $I_{\lambda,\beta}$ satisfies the (PS) condition.
\begin{lemma}\label{l4.6} For fixed $(\lambda,\beta)\in (0,1]\times(0,1]$, $I_{\lambda,\beta}$ satisfies the (PS)$_c$ condition with any $c\in\mathbb{R}$.
\end{lemma}
{\bf Proof}:  Assume that there exists $\{u_n\}\subset E$ satisfies
\begin{equation}\label{4.2.8}I_{\lambda,\beta}(u_n)\rightarrow c\ \text{and}\ I'_{\lambda,\beta}(u_n)\rightarrow 0,\end{equation}
 as $n\rightarrow\infty$. Set $\gamma\in (4,r)$, for any $\epsilon>0$, by (\ref{2.2}) there exists $C_\epsilon>0$ such that
 $$\aligned &I_{\lambda,\beta}(u_n)-\frac{1}{\gamma}\langle I'_{\lambda,\beta}(u_n),\frac{f(u_n)}{f'(u_n)}\rangle\\
=& \int_{\mathbb{R}^3}\Bigl(\frac12-\frac{1+4f^2(u_n)}{{\gamma}(1+2f^2(u_n))}
\Bigr)|\nabla {u_n}|^2_2+\Bigl(\frac12-\frac1{\gamma}\Bigr)\int_{\mathbb{R}^3}V(x)f^2(u_n)+
\frac{(r-\gamma)\beta}{r\gamma}|f(u_n)|^r_r\\&+\Bigl(\frac{1}{2(1+\alpha)}
-\frac{1}{\gamma}\Bigr)\lambda
|f(u_n)|^{2(1+\alpha)}_2
+\int_{\mathbb{R}^3}\Bigl[\frac1\gamma g(f(u_n))f(u_n)-G(f(u_n))\Bigr]
\\ \geq &\Bigl(\frac12-\frac2{\gamma}-\epsilon\Bigr)\bigl[|\nabla u_n|^2_2+\int_{\mathbb{R}^3}V(x)f^2(u_n)\bigr]+\Bigl(\frac{1}{2(1+\alpha)}-\frac{1}{\gamma}\Bigr)\lambda|f(u_n)|^{2(1+\alpha)}_2
\\&-
C_\epsilon|f(u_n)|^p_p+\frac{(r-\gamma)\beta}{r\gamma}|f(u_n)|^r_r.\endaligned$$
Choosing small enough $\epsilon>0$ we have
$$\aligned
&I_{\lambda,\beta}(u_n)-\frac{1}{\gamma}\langle I'_{\lambda,\beta}(u_n),\frac{f(u_n)}{f'(u_n)}\rangle\\
\geq&\Bigl(\frac14-\frac1{\gamma}\Bigr)\bigl[|\nabla u_n|^2_2+\int_{\mathbb{R}^3}V(x)f^2(u_n)\bigr]+\Bigl(\frac{1}{2(1+\alpha)}-\frac{1}{\gamma}\Bigr)\lambda|f(u_n)|^{2(1+\alpha)}_2
\\&-
C_1|f(u_n)|^p_p+\frac{(r-\gamma)\beta}{r\gamma}|f(u_n)|^r_r.
\endaligned$$
Using $\alpha<1$ and (\ref{4.1.5}) we infer

$$|\nabla u_n|^2_2+\int_{\mathbb{R}^3}V(x)f^2(u_n)\\
\leq  C\Bigl(C_2+|I_{{\lambda,\beta}}(u_n)|+o_n(1)\bigl\|\frac{f(u_n)}{f'(u_n)}\bigr\|\Bigr).$$
From (\ref{3.0.1}) and (\ref{2.3.4}) it follows that
\begin{equation*}|\nabla u_n|^2_2+\int_{\mathbb{R}^3}V(x)f^2(u_n)
\leq C\Bigl[C_3+o_n(1)\bigl(1+|\nabla u_n|_2+|f(u_n)|^2_2\bigr)\Bigr].\end{equation*}
Then $\{|\nabla u_n|^2_2+|f(u_n)|^2_2\}$ is bounded and so $\{u_n\}$ is bounded in $E$ using (\ref{2.3.4}) again.

Suppose $u_n\rightharpoonup u_0$ in $E$ and $u_n\rightarrow u_0$ in $L^s(\mathbb{R}^3)$ with $2<s<6$.
It suffices to show that $u_n\rightarrow u_0$ {in} $E$ up to a subsequence.
Indeed, arguing as (\ref{4.2.0}) and (\ref{4.2.1}) we get
 $$\aligned&|f(u_n)|^{2\alpha}_2\int_{\mathbb{R}^3}\frac{f(u_n)f'(u_n)}{u_n}u_n(u_n-u_0)
-|f(u_0)|^{2\alpha}_2\int_{\mathbb{R}^3}\frac{f(u_0)f'(u_0)}{u_0}u_0(u_n-u_0)
 \\=&|f(u_n)|^{2\alpha}_2\int_{\mathbb{R}^3}\frac{f(u_n)f'(u_n)}{u_n}(u_n-u_0)^2+o_n(1),\endaligned$$
 and $$\int_{\mathbb{R}^3}(\bar{g}(x,u_n)
-\bar{g}(x,u_0))(u_n-u_0)=o_n(1).$$Then
 $$\aligned &\langle I'_{\lambda,\beta}(u_n)-I'_{\lambda,\beta}(u_0),u_n-u_0\rangle\\
=&\|u_n-u_0\|^2+\lambda|f(u_n)|^{2\alpha}_2\int_{\mathbb{R}^3}\frac{f(u_n)f'(u_n)}{u_n}(u_n-u_0)^2+o_n(1).
\endaligned$$
On the other hand, by (\ref{4.2.8}) we obtain $\langle I'_{\lambda,\beta}(u_n)-I'_{\lambda,\beta}(u_0),u_n-u_0\rangle=o_n(1)$. Therefore, $u_n\rightarrow u_0$ in $E$.\ \ \ \ $\Box$

\subsection{Existence of a sign-changing solution}

In this subsection, we apply Theorem \ref{th1} to prove the existence of
sign-changing solutions of the problem (\ref{3.0.0}), and take $X=E$, $P=P^+_\epsilon$, $Q=Q^+_\epsilon$ and $J=I_{\lambda,\beta}$. We
will prove that $\{P^+_\epsilon, Q^+_\epsilon\}$ is an admissible family of invariant sets for the functional
 $I_{\lambda,\beta}$ at any level $c\in\mathbb{R}$. Indeed, if $K_c\backslash W=\emptyset$, then $K_c\subset W$. By Lemma \ref{l4.6}, the functional $I_{\lambda,\beta}$
satisfies the (PS)$_c$ condition, and so $K_c$ is compact. Then $2\delta:=dist(K_c,\partial W)>0$.

Here we give a deformation lemma to the functional $I_{\lambda,\beta}$ whose proof is almost the same
as that of \cite[Lemma 3.6]{LWZ}.

\begin{lemma} \label{l4.8}(Deformation lemma) If $K_c\backslash W=\emptyset$, then there exists $\epsilon_0>0$ such that, for $0<\epsilon<\epsilon'<\epsilon_0$, there
exists a continuous map $\sigma:[0,1]\times E\rightarrow E$ satisfying

\noindent(1) $\sigma(0,u)=u$ for $u\in E$;\\
\noindent(2) $\sigma(t,u)=u$ for $t\in[0,1]$, $u\not\in I^{-1}_{\lambda,\beta}[c-\epsilon',c+\epsilon']$;\\
\noindent(3) $\sigma(1,I_{\lambda,\beta}^{c+\epsilon}\backslash W)\subset I^{c-\epsilon}$;\\
\noindent(4) $\sigma(t,\overline{P^+_\epsilon})\subset \overline{P^+_\epsilon}$ and $\sigma(t,\overline{P^-_\epsilon})\subset \overline{P^-_\epsilon}$ for $t\in[0,1]$.
\end{lemma}
\begin{lemma}\label{l4.9}For any  $q\in[2,6]$, there exists $C>0$ independent of $\epsilon$ such that $|u|_q\leq C\epsilon$ for $u\in M=P^+_\epsilon\cap P^-_\epsilon$.\end{lemma}
{\bf Proof}: For any fixed $u\in M$, we have
$$|u^{\pm}|_q=\inf_{v\in P^{\mp}}|u-v|_q\leq C\inf_{v\in P^{\mp}}\|u-v\|\leq C dist(u,P^{\mp}).$$
Then $|u|_q\leq C\epsilon$ for $u\in M$.\ \ \ \ \ $\Box$

\begin{lemma}\label{l4.10}If $\epsilon>0$ small enough, then $I_{\lambda,\beta}(u)\geq\frac{\epsilon^2}{12}$ for $u\in \Sigma=\partial P^+_\epsilon\cap\partial P^-_\epsilon$, and so $c_*:=\inf_{u\in \Sigma } I_{\lambda,\beta}(u)\geq\frac{\epsilon^2}{12}$.\end{lemma}
{\bf Proof}: For $u\in \partial P^+_\epsilon\cap\partial P^-_\epsilon$, we have
$\|u^{\pm}\|\geq dist(u, P^{\pm})=\epsilon$. In view of (\ref{4.0.2}) and Lemma \ref{l4.9}, for any $\delta>0$, there exists $C_\delta>0$ such that
$$ I_{\lambda,\beta}(u)\geq \frac12\|u\|^2-\int_{\mathbb{R}^3}\bar{{G}}(x,u)\geq  \bigl(\frac12-\delta\bigr)\|u\|^2-C_\delta|u|^{\frac r2}_{\frac r2}.$$
Choosing $\delta$ small enough, we get
$$ I_{\lambda,\beta}(u)\geq\frac16\|u\|^2-C|u|^{\frac r2}_{\frac r2}
\geq\frac16\epsilon^2-C\epsilon^{\frac r2}\geq\frac{\epsilon^2}{12},$$
for $\epsilon$ small enough.\ \ \ \ $\Box$

{\bf Proof of Theorem 1.1 (Existence part)} In what follows, we divide the proof into three steps.

{\bf Step 1}. In this step, we will apply Theorem \ref{th1} to look for a sign-changing solutions of (\ref{3.0.0}) for any fixed $(\lambda,\beta)\in(0,1]\times(0,1]$. By Lemma \ref{l4.8}, we know $\{P^+_\epsilon, P^-_\epsilon\}$ is an admissible family of invariant sets with respect
to $I_{\lambda,\beta}$ at any level $c\in\mathbb{R}$. It suffices to verify assumptions (1)-(3) of
Theorem \ref{th1}.

Let $v_1,v_2\in C^\infty_0(\mathbb{R}^3)\backslash\{0\}$ be such that $\text{supp}(v_1)\cap \text{supp}(v_2)=\emptyset$ and $v_1\leq0$, $v_2\geq0$. For $(t,s)\in \triangle$, let
$$\varphi_0(t,s)(\cdot)=R(tv_1(R^{-1}\cdot)+sv_2(R^{-1}\cdot)),$$
where $R>0$ will be determined later. Obviously, for $t,s\in[0,1]$,
$\varphi_0(0,s)=Rsv_2(R^{-1}\cdot)\in P^+_\epsilon$ and $\varphi_0(t,0)=Rtv_1(R^{-1}\cdot)\in P^-_\epsilon$.
 Observe that $\rho:=\min\{|tv_1+(1-t)v_2|_2:0\leq t\leq1\}>0$. Then $|u|^2_2\geq \rho^2 R^5$ for $u\in \varphi_0(\partial_0\triangle)$. It follows from Lemma \ref{l4.9} that $\varphi_0(\partial_0\triangle)\cap M=\emptyset$.  In view of Lemma \ref{l4.10}, for small $\epsilon$ we have $c_*=\inf_{u\in \Sigma}I_{\lambda,\beta}(u)\geq\frac{\epsilon^2}{12}$ for any $(\lambda,\beta)\in(0,1]\times(0,1]$. Below we show that $\sup_{u\in \varphi_0(\partial_0\triangle)}I_{\lambda,\beta}(u)<0$. Set $u_t=\varphi_0(t,1-t)$ for $t\in[0,1]$.
A direct computation shows that
\begin{equation*}\aligned
 |\nabla u_t|^2_2&=R^3\int_{\mathbb{R}^3}(t^2|\nabla v_1|^2+(1-t)^2|\nabla v_2|^2);\\
|f(u_t)|^{2(1+\alpha)}_2&=R^{3(1+\alpha)}\Bigl(
\int_{\mathbb{R}^3}(|f(Rtv_1)|^2+|f(R(1-t)v_2)|^2)\Bigr)^{1+\alpha};\\
|f(u_t)|^2_2&=R^{3}
\int_{\mathbb{R}^3}(|f(Rtv_1)|^2+|f(R(1-t)v_2)|^2);\\
|f(u_t)|^\mu_\mu&=R^{3}
\int_{\mathbb{R}^3}(|f(Rtv_1)|^\mu+|f(R(1-t)v_2)|^\mu).
\endaligned\end{equation*}
Note that $G(s)\geq C|s|^\mu-C_1$ for any $s\in\mathbb{R}$, we deduce
$$\aligned I_{\lambda,\beta}(u_t)\leq&\frac{R^3}{2}\int_{\mathbb{R}^3}\bigl(t^2|\nabla v_1|^2+(1-t)^2|\nabla v_2|^2\bigr)+
\frac{\lambda R^{3(1+\alpha)}}{2(1+\alpha)}\bigl(|f(Rtv_1)|^{2(1+\alpha)}_2
\\&+|f(R(1-t)v_2)|^{2(1+\alpha)}_2\bigr)-\int_{\text{supp}v_1}G\bigl(f(Rtv_1)\bigr)-
\int_{\text{supp}v_2}G\bigl(f(R(1-t)v_2)\bigr)\\
\leq&\frac{R^3}{2}\int_{\mathbb{R}^3}\bigl(t^2|\nabla v_1|^2+(1-t)^2|\nabla v_2|^2\bigr)+
\frac{ R^{3(1+\alpha)}}{2(1+\alpha)}\bigl(|f(Rtv_1)|^{2(1+\alpha)}_2
\\&+|f(R(1-t)v_2)|^{2(1+\alpha)}_2\bigr)-
{R^3}\bigl(|f(Rtv_1)|^\mu_\mu+|f(R(1-t)v_2)|^\mu_\mu\bigr)+C_1R^3.\endaligned$$
Since
$$\frac{|f(s)|}{|s|^{\frac12}}\rightarrow 2^{\frac14},\ \text{as}\ |s|\rightarrow+\infty,\quad \text{and}\  \alpha<\frac{\mu-2}{3\mu+2}<\frac{1}{4}({\frac \mu2}-1),$$
 we have $I_{\lambda,\beta}(u_t)\rightarrow-\infty$ as $R\rightarrow+\infty$ uniformly for $(\lambda,\beta)\in(0,1]\times(0,1]$. So we can choose large enough $R>0$ independent on $\lambda$ and $\beta$ such that $I_{\lambda,\beta}(u_t)<0$. Then
$$\sup_{u\in \varphi_0(\partial_0\triangle)}I_{\lambda,\beta}(u)<0<\frac{\epsilon^2}{12}\leq c_*:=\inf_{u\in \Sigma}I_{\lambda,\beta}(u),\quad \forall (\lambda,\beta)\in (0,1]\times(0,1].$$
Applying Theorem \ref{th1} we know
$$c_{\lambda,\beta}=\inf_{\varphi\in \Gamma}\sup_{u\in \varphi(\triangle)\backslash W}I_{\lambda,\beta}(u),$$
 is a  critical value of $I_{\lambda,\beta}$ and \begin{equation}\label{4.11.0}c_{\lambda,\beta}\geq c_*\geq\frac{\epsilon^2}{12}.\end{equation} Therefore, there exists $u_{\lambda,\beta}\in E\backslash{(P^+_\epsilon \cup P^-_\epsilon)}$ such that $I_{\lambda,\beta}(u_{\lambda,\beta})=c_{\lambda,\beta}$ and $I'_{\lambda,\beta}(u_{\lambda,\beta})=0$ for fixed $(\lambda,\beta)\in (0,1]\times(0,1]$.

 {\bf Step 2}.  Set $\lambda\rightarrow0$ and $\beta\rightarrow0$. In view of the definition of $c_{\lambda,\beta}$,
 for any $(\lambda,\beta)\in(0,1]\times(0,1]$, there holds
\begin{equation}\label{4.1.9}c_{\lambda,\beta}\leq C_R:=\sup_{u\in \varphi_0(\triangle)}I_{1,0}(u)<+\infty,\end{equation}
where $C_R$ is  independent on $(\lambda,\beta)\in(0,1]\times(0,1]$.  Without loss of generality, we set $\lambda=\beta$. Choosing a sequence
$\{\lambda_n\}\subset(0,1]$ satisfying $\lambda_n\rightarrow 0^+$, we find a sequence
of sign-changing critical points $\{u_{\lambda_n}\}$ (still denoted by $\{u_{n}\}$ for simplicity) of $I_{\lambda_n,\beta_n}$ and $I_{\lambda_n,\beta_n}(u_n)=c_{\lambda_n,\beta_n}$. Now we show that $\{u_n\}$ is bounded in $E$. Note that
\begin{equation}\label{4.1.10}c_{\lambda_n,\beta_n}=\frac12|\nabla u_n|^2_2+\frac12 \int_{\mathbb{R}^3}V(x)f^2(u_n)+\frac{\lambda_n}{2(1+\alpha)}|f(u_n)|^{2(1+\alpha)}_2
-\int_{\mathbb{R}^3}G(f(u_n))-\frac{\beta_n}{r}|f(u_n)|^r_r,\end{equation}
and \begin{equation}\label{4.1.11}\aligned0=\langle I'_{\lambda_n,\beta_n}(u_n),\frac{f(u_n)}{f'(u_n)}\rangle=&\int_{\mathbb{R}^3}
\Bigl(1+\frac{2f^2(u_n)}{1+2f^2(u_n)}\Bigr)|\nabla u_n|^2+ \int_{\mathbb{R}^3}V(x)f^2(u_n)\\
&+{\lambda_n}|f(u_n)|^{2(1+\alpha)}_2-\int_{\mathbb{R}^3}g(f(u_n))f(u_n)
-{\beta_n}|f(u_n)|^r_r.\endaligned\end{equation}
Moreover, from Lemma \ref{l4.1.0} we have
\begin{equation}\label{4.1.12}\aligned &\frac{1}{2}|\nabla u_n|^2_2+\frac32\int_{\mathbb{R}^3}V(x)f^2(u_n)
+\frac12\int_{\mathbb{R}^3}(\nabla V(x),x) f^2(u_n)+\frac{3\lambda_n}{2(1+\alpha)}|f(u_n)|^{2(1+\alpha)}_2\\&-3\int_{\mathbb{R}^3}G(f(u_n))
-\frac{3\beta_n}{r}|f(u_n)|^r_r=0.\endaligned\end{equation}
Multiplying (\ref{4.1.10}), (\ref{4.1.11}) and (\ref{4.1.12}) by $4$, $-\frac{1}{\mu}$ and $-1$ respectively and adding them up,
we obtain
$$\aligned 4c_{\lambda_n,\beta_n}=&\int_{\mathbb{R}^3}
\Bigl[\frac32-\frac1\mu\bigl(1+\frac{2f^2(u_n)}{1+2f^2(u_n)}\bigr)
\Bigr]|\nabla u_n|^2+\bigl(\frac12-\frac1\mu\bigr)
\int_{\mathbb{R}^3}V(x)f^2(u_n)\\&-\frac12\int_{\mathbb{R}^3}(\nabla V(x),x) f^2(u_n)+\bigl[\frac1{2(1+\alpha)}-\frac1\mu\bigr]
\lambda_n|f(u_n)|^{2(1+\alpha)}_2\\
&+\int_{\mathbb{R}^3}\bigl[-G(f(u_n))+\frac1\mu g(f(u_n))f(u_n)\bigr]+\bigl(\frac1\mu-\frac1r\bigr)\beta_n|f(u_n)|^r_r.
\endaligned$$
Since $\alpha<\frac{\mu-2}{3\mu+2}<\frac{\mu}{2}-1$, it follows from (V$_2$),
($g_3$) and (\ref{4.1.9}) that
$$4C_R\geq\bigl(\frac32-\frac2\mu\bigr)|\nabla u_n|^2_2,$$
which implies that there exists $C>0$ independent of $\lambda,\beta$ such that
\begin{equation}\label{4.1.13}|\nabla u_n|^2_2<C.\end{equation}
Moreover, in view of (\ref{4.1.9}) and (\ref{4.1.10})
we infer that for small $\delta>0$, there exists $C_\delta>0$ such that
\begin{equation}\label{4.1.14}\aligned C_R&\geq\frac12\int_{\mathbb{R}^3}V(x)f^2(u_n)
-\int_{\mathbb{R}^3}G(f(u_n))-\frac{\beta_n}{r}|f(u_n)|^r_r\\
&\geq\frac{1-\delta}2\int_{\mathbb{R}^3}V(x)f^2(u_n)
-C_\delta\int_{\mathbb{R}^3}|u_n|^{6}-\frac{1}{r}|f(u_n)|^r_r.\endaligned\end{equation}
Observe that $r\in (4,12)$, we may assume that $r=2\tau+12(1-\tau)$, $\tau\in (0,1)$. Then for the above $\delta$, there exists $\bar{C}_\delta>0$ such that
\begin{equation}\label{4.1.15}\aligned |f(u_n)|^r_r&\leq |f(u_n)|^{2\tau}_2|f(u_n)|^{12(1-\tau)}_{12}\leq |f(u_n)|^{2\tau}_2|u_n|^{6(1-\tau)}_6\\
&\leq \delta \int_{\mathbb{R}^3}V(x)f^2(u_n)+\bar{C}_\delta|u_n|^{6}_6.\endaligned\end{equation}
According to (\ref{4.1.14}) and (\ref{4.1.15}) we have
$$C_R\geq\frac{1-2\delta}2\int_{\mathbb{R}^3}V(x)f^2(u_n)
-(C_\delta+\bar{C}_\delta)\int_{\mathbb{R}^3}|u_n|^{6}.$$
Using (\ref{4.1.13}) and the arbitrariness  of $\delta$, we get $|f(u_n)|^2_2$ is bounded. By (\ref{2.3.4}) we know $\{u_n\}$ is bounded in $E$. From  (\ref{4.11.0}) it follows that
\begin{equation}\label{4.18.0}\aligned\lim_{n\rightarrow\infty}I(u_n)&=\lim_{n\rightarrow\infty}\Bigl(
I_{\lambda_n,\beta_n}(u_n)
-\frac{\lambda_n}{2(1+\alpha)}|f(u_n)|^{2(1+\alpha)}_2
+\frac{\beta_n}{r}|f(u_n)|^r_r\Bigr)\\&=\lim_{n\rightarrow\infty}
c_{\lambda_n,\beta_n}:=c^*\geq c_*\geq\frac{\epsilon^2}{12}.\endaligned\end{equation}
Moreover, for any $\psi\in C^\infty_0(\mathbb{R}^3)$,
$$\aligned\lim_{n\rightarrow\infty}\langle I'(u_n),\psi\rangle=&\lim_{n\rightarrow\infty}\Bigl(\langle I'_{\lambda_n,\beta_n}(u_n),\psi\rangle-
\lambda_n|f(u_n)|^{2\alpha}_2\int_{\mathbb{R}^3}f(u_n)f'(u_n)\psi\\&+
\beta_n\int_{\mathbb{R}^3}|f(u_n)|^{r-2}f(u_n)f'(u_n)\psi\Bigr)=0.
\endaligned$$
Then $\{u_n\}$ is a bounded (PS) sequence for $I$ at level $c^*$. Hence, we may assume that $u_n\rightharpoonup u^*$ in $E$ and $u_n\rightarrow u^*$ in $L^s(\mathbb{R}^3)$ with $s\in (2,6)$. By Lemma \ref{l3.1.0} (3), arguing as (\ref{4.2.1}) we have
 $$\int_{\mathbb{R}^3}(\tilde{g}(x,u_n)
-\tilde{g}(x,u^*))(u_n-u^*)=o_n(1).$$
 Then
 $$\langle I'(u_n)-I'(u^*),u_n-u^*\rangle
=\|u_n-u^*\|^2+o_n(1).$$
 On the other hand, since $I'(u_n)\rightarrow0$ we get $\langle I'(u_n)-I'(u^*),u_n-u^*\rangle=o_n(1)$. Consequently
\begin{equation}\label{4.2.4}u_n\rightarrow u^*\ \text{in}\ E.\end{equation}
Then
 $I'(u^*)=0$ and $I(u^*)=c^*$. From the fact that $u_n\in E\backslash{(P^+_\epsilon\cup P^-_\epsilon)}$ we infer $u^*\in E\backslash{(P^+_\epsilon\cup P^-_\epsilon)}$ and so $u^*$
is a sign-changing solution of problem (\ref{2.0.0}).

{\bf Step 3}. Define
$$\bar{c}:=\inf_{u\in \mathcal{N}}I(u), \ \ \mathcal{N}:=\{u\in E: I'(u)=0, u^{\pm}\not\equiv0\}.$$
Based on Step 2, we see that $\mathcal{N}\neq\emptyset$ and $\bar{c}\leq c^*$, where $c^*$ is given in (\ref{4.18.0}). By the definition of $\bar{c}$, there exists $\{u_n\}\subset E$ such that $I(u_n)\rightarrow \bar{c}$ and $I'(u_n)=0$. Arguing as (\ref{4.2.4}),  we infer that there exists $u\neq0$ such that $u_n\rightarrow u$ in $E$, $I(u)=\bar{c}$ and $I'(u)=0$. Moreover, note that $\langle I'(u_n), u^{\pm}_n\rangle=0$.
From Lemma \ref{l3.1.0} (3) we infer that for any $\delta>0$, there exists $C_\delta>0$ such that
$$ |\nabla u^{\pm}_n|^2_2+\int_{\mathbb{R}^3}V(x)(u^{\pm}_n)^2
=\int_{\mathbb{R}^3}\tilde{g}(x,u^{\pm}_n)u^{\pm}_n
\leq \delta\|u^{\pm}_n\|^2+C_\delta\|u^{\pm}_n\|^q.$$
Then
$\|u^{\pm}_n\|\geq \varrho>0$ and $\|u^{\pm}\|\geq \frac{\varrho}{2}$.
Consequently, $u$ is a least energy sign-changing solution of problem (\ref{2.0.0}). The proof is complete.\ \ \ \ $\Box$

\subsection{Multiplicity of sign-changing solutions}
 In this subsection, $g$ is assumed to be odd, and so $I_{\lambda,\beta}$ is even. We shall apply Theorem \ref{th2} to obtain infinitely many sign-changing solutions and set $X=E$, $G=-id$, $J=I_{\lambda,\beta}$ and $P=P^+_\epsilon$. Then $M=P^+_\epsilon\cap P^-_\epsilon$, $\Sigma=\partial P^+_\epsilon\cap \partial P^-_\epsilon$, and $W=P^+_\epsilon\cup P^-_\epsilon$. Since $K_c$ is compact, there exists a
symmetric open neighborhood $N$ of $K_c\backslash W$
such that $\gamma(\overline{N})<+\infty$.

\begin{lemma}\label{l4.11} There exists $\epsilon_0>0$ such that, for $0<\epsilon<\epsilon'<\epsilon_0$, there
exists a continuous map $\sigma:[0,1]\times E\rightarrow E$ satisfying\\
\noindent(1) $\sigma(0,u)=u$ for $u\in E$;\\
\noindent(2) $\sigma(t,u)=u$ for $t\in[0,1]$, $u\not\in I^{-1}[c-\epsilon',c+\epsilon']$;\\
\noindent(3) $\sigma(1,I^{c+\epsilon}\backslash (N\cup W))\subset I^{c-\epsilon}$;\\
\noindent(4) $\sigma(t,\overline{P^+_\epsilon})\subset \overline{P^+_\epsilon}$ and $\sigma(t,\overline{P^-_\epsilon})\subset \overline{P^-_\epsilon}$ for $t\in[0,1]$;\\
\noindent(5) $\sigma(t,-u)=-\sigma(t,u)$ for $(t,u)\in[0,1]\times E$.\end{lemma}
{\bf Proof}: The arguments of (1)-(4) are similar to those of Lemma \ref{l4.8}. Concerning (5), since $I_{\lambda,\beta}$ is even and the operator $B_{\lambda,\beta}$ in Lemma \ref{l4.5} is odd,
we infer $\sigma$ is odd in $u$.\ \ \ \ $\Box$

{\bf Proof of Theorem 1.1} ({\bf Multiplicity part}) We complete the proof in two steps.

{\bf Step 1}. Since $g$ is odd, in view of Lemma \ref{l4.11} we know
$P^+_\epsilon$ is a $G-$admissible invariant set with respect to $I_{\lambda,\beta}$ at
level $c$. In order to apply Theorem \ref{th2}, we are now constructing $\varphi_n$. For any $n\in\mathbb{N}$, let $\{v_i\}^n_{i=1}\subset C^\infty_0(\mathbb{R}^3)\backslash\{0\}$ be such that $\text{supp}(v_i)\cap \text{supp}(v_j)=\emptyset$ for $i\neq j$. Define $\varphi_n\in C(B_n, E)$
as
$$\varphi_n(t)=R_n\Sigma^{n}_{i=1}t_iv_i(R^{-1}_n\cdot),\quad\ t=(t_1,t_2,...,t_n)\in B_n,$$
where $R_n>0$ will be determined later. Obviously, $\varphi_n(0)=0\in P^+_\epsilon\cap P^-_\epsilon$ and $\varphi_n(-t)=-\varphi_n(t)$ for $t\in B_n$.
Observe that
$$\rho_n=\min\{|t_1v_1+t_2v_2+...+t_nv_n|_2:\Sigma^{n}_{i=1}t^2_i=1\}>0,$$
then $|u_t|^2_2\geq \rho^2_n R^5_n$ for $u\in \varphi_n(\partial B_n)$ and it follows from Lemma \ref{l4.9} that $\varphi_n(\partial B_n)\cap (P^+_\epsilon\cap P^-_\epsilon)=\emptyset$. Similar to the proof of Theorem 1.1 (existence part), for large enough $R_n>0$ independent on $\lambda$ and $\beta$ we also have
$$\sup_{u\in \varphi_n(\partial B_n)}I_{\lambda,\beta}(u)<0<\inf_{u\in \Sigma}I_{\lambda,\beta}(u).$$ For any $j\in\mathbb{N}$ and $(\lambda,\beta)\in (0,1]\times(0,1]$, we define
$$c^j_{\lambda,\beta}=\inf_{B\in \Gamma_j}\sup_{u\in B\backslash W}I_{\lambda,\beta}(u),$$
where $W=P^+_\epsilon\cup P^-_\epsilon$ and $\Gamma_j$ is defined in Theorem \ref{th2}. Applying Theorem \ref{th2} and Lemma \ref{l4.10},
for any $(\lambda,\beta)\in (0,1]\times(0,1]$ and $j\geq2$,
\begin{equation}\label{4.3.2}\frac{\epsilon^2}{12}\leq\inf_{u\in \Sigma}I_{\lambda,\beta}(u)\leq c^j_{\lambda,\beta}\rightarrow+\infty, \ \ \text{as} \ j\rightarrow+\infty,\end{equation}
and there exists $\{u^j_{\lambda,\beta}\}\subset E\backslash W$ such that $I_{\lambda,\beta}(u^j_{\lambda,\beta})=c^j_{\lambda,\beta}$ and $I'_{\lambda,\beta}(u^j_{\lambda,\beta})=0$.

{\bf Step 2}. Using similar arguments as those in Theorem 1.1 (existence part), for any fixed $j\geq2$, $\{u^j_{\lambda,\beta}\}_{\lambda,\beta\in(0,1]}$ is bounded in $E$. Namely, there exists $C>0$ independent of $\lambda,\beta$ such that $\|u^j_{\lambda,\beta}\|\leq C$. Without loss of generality, we assume $u^j_{\lambda,\beta}\rightharpoonup u^j_*$ in $E$ as $\lambda,\beta\rightarrow0$.  By (\ref{4.3.2}) we have
$$\frac{\epsilon^2}{12}\leq  c^{j}_{\lambda,\beta}\leq c_{R_n}:=\sup_{u\in \phi_n(B_n)}I_{1,0}(u),$$
where $c_{R_n}$ is independent of $\lambda,\beta$. Assume that $c^{j}_{\lambda,\beta}\rightarrow c^j_*$ as $\lambda,\beta\rightarrow0$. Then as the proof of Theorem 1.1 (existence part), we can prove that $u^j_{\lambda,\beta}\rightarrow u^j_*$ in $E$ as $\lambda,\beta\rightarrow0^+$ and $u^j_*\in E\backslash W$ satisfying $I'(u^j_*)=0$ and $I(u^j_*)=c^j_*$. Below we claim that $c^j_*\rightarrow+\infty$ as $j\rightarrow\infty$. Indeed, by ($g_1$) and ($g_2$) we infer
\begin{equation}\label{4.1.16}\aligned I_{\lambda,\beta}(u)&\geq\frac12|\nabla u|^2_2+\frac12\int_{\mathbb{R}^3}V(x)f^2(u)
-\int_{\mathbb{R}^3}G(f(u))-\frac1r|f(u)|^r_r\\
&\geq\frac12|\nabla u|^2_2+\frac12\int_{{\mathbb{R}}^3}V(x)f^2(u)
-\int_{\mathbb{R}^3}\bigl(\frac{\inf_{\mathbb{R}^3}V}4f^2(u)+\frac{C}{r}|f(u)|^r\bigr)
-\frac1r|f(u)|^r_r\\
&\geq\frac12|\nabla u|^2_2+\frac14\int_{\mathbb{R}^3}V(x)f^2(u)
-\frac{1+C}{r}|f(u)|^r_r:=\bar{I}(u),
\endaligned\end{equation}
where $C>0$ depends on $\inf_{\mathbb{R}^3}V$. Since $r\in(4,12)$, it is easy to verify that the (PS) condition of $\bar{I}$ holds. Then replacing $I_{\lambda,\beta}$ by $\bar{I}$,  the arguments of $I_{\lambda,\beta}$ are still valid with some suitable modifications. In particular, the assumptions of Theorem \ref{th2} are satisfied. So
we can define
\begin{equation*} d^j:=\inf_{B\in \Gamma_j}\sup_{u\in B\backslash W}\bar{I}(u),\end{equation*}
where $W=P^+_\epsilon\cup P^-_\epsilon$ and $\Gamma_j$ is defined in
 Theorem \ref{th2}. Moreover, Theorem \ref{th2} implies that
$d^j\rightarrow+\infty$ as $j\rightarrow+\infty$. Combining (\ref{4.1.16}) and the definitions of $c^{j}_{\lambda,\beta}$ and $d^j$, we have $c^{j}_{\lambda,\beta}\geq d^j$. Taking $\lambda,\beta\rightarrow0^+$ we get $c^j_*\geq d^j\rightarrow+\infty$ as $j\rightarrow+\infty$. Hence, problem (\ref{2.0.0}) has infinitely many sign-changing solutions. This
ends the proof. \ \ \ \ \ $\Box$

\section{Proof of Theorem 1.2}
\renewcommand{\theequation}{4.\arabic{equation}}
In this section, we are devoted to showing Theorem 1.2 and assume (V$_1$), ($g_1$), ($g_2$), ($g_4$) and ($g_5$) are satisfied. As said in Remark 1.3, we shall prove Theorem 1.2 without using the perturbation  method. However, under the assumptions of Theorem 1.2, we will need to show that the functional satisfies the (Ce) condition and establish the deformation lemma under the (Ce) condition.

\subsection{Properties of operator $A$}
We firstly introduce an auxiliary operator $A$, which will be used to construct the descending flow for
the functional $I$ given in (\ref{2.3}). As  an application of Lax-Milgram theorem, for any $u\in E$, there is a unique solution $v=A(u)\in E$
of the equation
\begin{equation}\label{3.1} -\Delta v+V(x)v=\tilde{g}(x,u),\end{equation}
where $\tilde{g}(x,u)$ is given in (\ref{2.1}).
Arguing as Lemma \ref{l4.2}, we have that $A$ is continuous and compact.
\begin{lemma}\label{l3.2}(1) $\langle I'(u), u-A(u)\rangle\geq\|u-A(u)\|^2$ for all $u\in E$;\\
(2) $\|I'(u)\|\leq\|u-A(u)\|$ for all $u\in E$;\\
(3) for $a<b$ and $\alpha>0$, there exists $\beta>0$ such that $\|u-A(u)\|\geq\beta$ if $u\in E$, $I(u)\in [a,b]$ and $\|I'(u)\|\geq\alpha$.\end{lemma}
{\bf Proof}: (1) For any $u\in E$, we have $$\langle I'(u), u-A(u)\rangle=\langle I'(u)-I'(A(u)), u-A(u)\rangle=\|u-A(u)\|^2.$$

\noindent(2) For any $u, \varphi\in E$, we get
$$\aligned\langle I'(u), \varphi\rangle=\langle I'(u)-I'(A(u)), \varphi\rangle=\langle u-A(u), \varphi\rangle\leq\|u-A(u)\|\|\varphi\|.\endaligned$$
Then $\|I'(u)\|\leq\|u-A(u)\|$.

\noindent(3) From (2) it follows that the conclusion (3) holds true.\ \ \ \ \ $\Box$

\subsection{Invariant subsets of descending flows}
The notations $P^+, P^-, P^+_\epsilon, P^-_\epsilon$ and $W$ in Section 3.2 are still valid. By Lemma \ref{l3.1.0} (3), arguing as Lemma \ref{l4.4} with small modifications, we have:

\begin{lemma}\label{l3.4} There is $\epsilon_0>0$ such that, for $\epsilon\in (0,\epsilon_0)$,\\
\noindent(1) $A(\partial P^{-}_\epsilon)\subset P^{-}_\epsilon$ and every nontrivial solution $u\in P^{-}_\epsilon$ is negative.

\noindent(2) $A(\partial P^{+}_\epsilon)\subset P^{+}_\epsilon$ and every nontrivial solution $u\in P^{+}_\epsilon$ is positive.
\end{lemma}

Denote the set of fixed points of $A$ by $K$. Similar to Lemma \ref{l4.5}, we have the following result.
\begin{lemma}\label{l3.5}There exists a locally Lipschitz continuous operator $B: E_0=E\backslash K\rightarrow E$ such
that \\
(1)\ $B(\partial P^+_\epsilon)\subset P^+_\epsilon$ and $B(\partial P^-_\epsilon)\subset P^-_\epsilon$ for $\epsilon\in (0,\epsilon_0)$;

\noindent(2)\ $\frac12\|u-B(u)\|\leq\|u-A(u)\|\leq 2\|u-B(u)\|$ for all $u\in E_0$;

\noindent(3)\ $\langle I'(u), u-B(u)\rangle\geq\frac12\|u-A(u)\|^2$ for all $u\in E_0$;

\noindent(4)\ if $g$ is odd, then $B$ is odd.
\end{lemma}

 Below we show that, $I$ satisfies the (Ce) condition. Although the proof is inspired by \cite{Jeanjean}, some crucial modifications are needed to overcome the difficulties caused by the change $f$.
\begin{lemma}\label{l5.1} For any $c\in\mathbb{R}$, $I$ satisfies the (Ce)$_c$ condition.
\end{lemma}
{\bf Proof}: Let $\{u_n\}\subset E$ be a (Ce)$_c$ sequence of $I$, i.e.
\begin{equation}\label{5.1}I(u_n)\rightarrow c\ \text{and}\ \|u_n\|(1+\|I'(u_n)\|)\rightarrow0.\end{equation}
Firstly, we show that $\|u_n\|$ is bounded. Argue by contradiction we may assume that $\|u_n\|\rightarrow+\infty$. Setting $v_n:=\frac{u_n}{\|u_n\|}$. Up to a subsequence, we suppose that
$v_n\rightharpoonup v$ in $E$, $v_n\rightarrow v$ in $L^s(\mathbb{R}^3)$ with $2<s<6$ and $v_n(x)\rightarrow v(x)$ a.e. in $\mathbb{R}^3$.
Below we consider two cases that $v=0$ and $v\neq0$ separately, and show that there will be a  contradiction in both cases.

{\bf Case 1}: $v=0$.

In this case, let $t_n\in[0,1]$ such that
$ I(t_nu_n)=\max_{t\in[0,1]}I(tu_n).$
For any positive constant $M$, we have $2M^{\frac12}{\|u_n\|^{-1}}\in(0,1)$ for large $n$. Denote $\tilde{v}_n=2M^{\frac12}v_n$.
From Lemma \ref{l3.1.0} (3) and Lebesgue dominated convergence theorem it follows that
$\int_{\mathbb{R}^3}\tilde{G}(x,\tilde{v}_n)=o_n(1).$
Then for $n$ large we get
$$I(t_nu_n)\geq I(\tilde{v}_n)=\frac{1}{2}\|\tilde{v}_n\|^2
-\int_{\mathbb{R}^3}\tilde{G}(x,\tilde{v}_n)\geq M.$$
Hence
\begin{equation}\label{5.3}\lim_{n\rightarrow\infty}I(t_nu_n)=+\infty.\end{equation}
On the other hand, note that $I(0)=0$ and $I(u_n)\rightarrow c$, we have $t_n\in (0,1)$ and then $\frac{d}{dt}I(tu_n)|_{t=t_n}=0$.
Denote $$\bar{f}(t):=f^2(t)-f(t)f'(t)t,\quad\ \bar{\bar{f}}(t):=f(t)f'(t)t-\frac{f^2(t) }2.$$
Clearly,  using Lemma \ref{l2.1} (6) we have $\bar{f}(t), \bar{\bar{f}}(t)\geq0$. We claim that
\begin{equation}\label{5.4}\bar{f}(t)\ \text{and} \ \bar{\bar{f}}(t)\text{\ are nondecreasing in}\ (0,+\infty)\ \text{and nonincreasing in}\ (-\infty,0).\end{equation}
Indeed, by Lemma \ref{l2.1} (6) and (8), for any $t>0$ we deduce $$\bar{f}'(t)=f(t)f'(t)-[f'(t)]^2t+2f^2(t)[f'(t)]^4t\geq2f^2(t)[f'(t)]^4t\geq0,$$
and
$$\bar{\bar{f}}'(t)=t[f'(t)]^2+f(t)f''(t)t=t[f'(t)]^2[1-2f^2(t)(f'(t))^2]\geq0.$$
Note that $\bar{f}$ and $\bar{\bar{f}}$ are even, so (\ref{5.4}) holds true. Using (\ref{5.4}) and ($g_5$), we infer
\begin{align}\label{4.2.5}
I(t_nu_n)=&I(t_nu_n)-\frac{1}{2}\langle I'(t_nu_n),t_nu_n\rangle\nonumber\\
=&\frac12\int_{\mathbb{R}^3}V(x)\bigl[f^2(t_nu_n)-f(t_nu_n)f'(t_nu_n)t_nu_n\bigr]\nonumber\\
&+\int_{\mathbb{R}^3}\Bigl[\frac12g(f(t_nu_n))f'(t_nu_n)t_nu_n-G(f(t_nu_n))\Bigr]\nonumber\\
=&\frac12\int_{\mathbb{R}^3}V(x)\bar{f}(t_nu_n)
+2\int_{\mathbb{R}^3}\frac{G(f(t_nu_n))}{f^2(t_nu_n)}\bar{\bar{f}}(t_nu_n)\nonumber\\
&+\int_{\mathbb{R}^3}\Bigl[\frac12g(f(t_nu_n))f(t_nu_n)-2G(f(t_nu_n))\Bigr]\frac{f'(t_nu_n)t_nu_n}{f(t_nu_n)}.
\end{align}
By Remark 1.2, $\frac{G(s)}{s^2}$ is nondecreasing for $s>0$ and nonincreasing for $s<0$. It follows that $\frac{G(f(t))}{f^2(t)}\le \frac{G(f(s))}{f^2(s)}$ if $|t|\le|s|$ and then
$$
\int_{\mathbb{R}^3}\frac{G(f(t_nu_n))}{f^2(t_nu_n)}\bar{\bar{f}}(t_nu_n)\le\int_{\mathbb{R}^3}\frac{G(f(u_n))}{f^2(u_n)}\bar{\bar{f}}(u_n).
$$
By Lemma \ref{l2.1}, $\frac{f'(t)t}{f(t)}\ge\frac12$ and $\frac{f'(t)t}{f(t)}\le1$ for any $t\not=0$. This yields that
$$
\frac{f'(t)t}{f(t)}\le 2\frac{f'(s)s}{f(s)},\mbox{for any}\,\,t,s\not=0.
$$
So we have
$$
\frac{f'(t_nu_n)t_nu_n}{f(t_nu_n)}\le 2\frac{f'(u_n)u_n}{f(u_n)},\,\, a. e.\,\, \mbox{in}\,\,\mathbb{R}^3.$$
Thanks to $(g_5)$ and the monotonicity of $\bar{f}$, by (\ref{4.2.5}) we get
\begin{align*}
I(t_nu_n)\le&\frac12\int_{\mathbb{R}^3}V(x)\bar{f}(u_n)
+2\int_{\mathbb{R}^3}\frac{G(f(u_n))}{f^2(u_n)}\bar{\bar{f}}(u_n)\nonumber\\
&+2\gamma\int_{\mathbb{R}^3}\Bigl[\frac12g(f(u_n))f(u_n)-2G(f(u_n))\Bigr]\frac{f'(u_n)u_n}{f(u_n)}\\
\le& 2\gamma[I(u_n)-\frac{1}{2}\langle I'(u_n),u_n\rangle].
\end{align*}
In view of (\ref{5.1}) we know
$I(t_nu_n)\leq 2\gamma[c+o_n(1)],$
contradicting (\ref{5.3}).

{\bf Case 2}: $v\not\equiv0$.

Set $\Omega=\{x\in \mathbb{R}^3: v(x)\neq0\}$. Then $|u_n(x)|\rightarrow\infty$ for $x\in \Omega$. From ($g_4$), Lemma \ref{l2.1} (5) and Fatou lemma it follows that
$$\aligned\liminf_{n\rightarrow\infty}\int_{\mathbb{R}^3}
\frac{G(f(u_n))}{\|u_n\|^2}&\geq
\liminf_{n\rightarrow\infty}\int_{\Omega}\frac{G(f(u_n))}{u^2_n}{v^2_n}\\&
=\liminf_{n\rightarrow\infty}\int_{\Omega}\frac{G(f(u_n))}{f^4(u_n)}
\frac{f^4(u_n)}{u^2_n}{v^2_n}=+\infty\endaligned$$
Hence
$$\aligned \frac{I(u_n)}{\|u_n\|^2}\leq\frac12-\int_{\mathbb{R}^3}
\frac{G(f(u_n))}{\|u_n\|^2}\rightarrow-\infty.\endaligned$$
However, by (\ref{5.1}) and $\|u_n\|\rightarrow+\infty$, we know $\frac{I(u_n)}{\|u_n\|^2}\rightarrow0$. This is a contradiction.

According to the above discussion, $\{u_n\}$ is bounded in $E$ and we assume that $u_n\rightharpoonup u$ in $E$. Then as the argument of (\ref{4.2.4}) we obtain $u_n\rightarrow u$ in $E$ after passing to a subsequence. This ends the proof.
\ \ \ \ $\Box$

\subsection{Existence of a sign-changing solution}

To apply Theorem \ref{th1}, we take $X=E$, $P=P^+_\epsilon$, $Q=Q^+_\epsilon$ and $J=I$. In the process of applying Theorem \ref{th1},  a deformation lemma of the functional $I$  corresponding to \cite[Lemma 3.6]{LWZ} is crucial. However, the proof of \cite[Lemma 3.6]{LWZ} depends on the fact that the functional satisfies the (PS) condition and here the functional $I$ merely satisfies the (Ce) condition, so we have to provide a new proof.

In the following, we firstly recall the local solvability and global existence theorems of solutions for initial value problems in \cite{Zhongchengkui}, see also \cite{Ligongbao}. Let $(X_0,\|\cdot\|_{X_0})$ be a real Banach space, suppose that $\varphi:[0,+\infty)\times X_0\rightarrow X_0$ is continuous and for all $x_0\in X_0$ and all $\alpha>0$, there exist $R>0$ and $L=L(x_0,\alpha,R)>0$ such that
\begin{equation}\label{5.6}\|\varphi(t,x)-\varphi(t,y)\|_{X_0}\leq L\|x-y\|_{X_0},\ \forall x,y\in B(x_0,R), t\in[0,\alpha]. \end{equation}
Consider the following initial value problem
\begin{equation}\label{5.7}\aligned
\left\{ \begin{array}{lll}
\frac{dx}{dt}=\varphi(t,x),\\
x(0)=x_0\in X_0,
\end{array}\right.\endaligned
\end{equation}

\begin{lemma}\label{l5.2}(\cite[Theorem 5.1]{Zhongchengkui})
Suppose that $\varphi$ satisfies the assumption (\ref{5.6}).
Then there exists $\beta>0$ such that the problem (\ref{5.7}) has a unique solution $x(t)$ in $[0,\beta]$ which
continuously depends on $x_0$. More generally, if $\|\varphi(x,t)-\varphi(y,t)\|_{X_0}\leq L\|x-y\|_{X_0}$, then
$$\|x(t)-y(t)\|_{X_0}\leq Le^{Lt}\|x_0-y_0\|_{X_0},\ \forall x,y\in X_0, t\in[0,\beta],$$
where $x(t)$ and $y(t)$ are the solutions of (\ref{5.7}) with initial values $x_0$ and $y_0$ respectively.
\end{lemma}

\begin{lemma}\label{l5.3}(\cite[Theorem 5.3]{Zhongchengkui})
Suppose that $\varphi$ satisfies the assumption (\ref{5.6}). If there exist $a,b>0$ such that
$$\|\varphi(t,x)\|_{X_0}\leq a+b\|x\|_{X_0}, \forall (t,x)\in[0,+\infty)\times {X_0}\mapsto {X_0},$$
then the unique local solution of (\ref{5.7}) can be extended as a global solution for $t\in[0,+\infty)$.
\end{lemma}

Now we are ready to state the deformation lemma under the (Ce) condition and give the proof.
\begin{lemma} \label{l5.4}(Deformation lemma) If $K_c\backslash W=\emptyset$, then there exists $\epsilon_0>0$ such that, for $0<\epsilon<\epsilon'<\epsilon_0$, there
exists a continuous map $\sigma:[0,1]\times E\rightarrow E$ satisfying

\noindent(1) $\sigma(0,u)=u$ for $u\in E$;\\
\noindent(2) $\sigma(t,u)=u$ for $t\in[0,1]$, $u\not\in I^{-1}[c-\epsilon',c+\epsilon']$;\\
\noindent(3) $\sigma(1,I^{c+\epsilon}\backslash W)\subset I^{c-\epsilon}$;\\
\noindent(4) $\sigma(t,\overline{P^+_\epsilon})\subset \overline{P^+_\epsilon}$ and $\sigma(t,\overline{P^-_\epsilon})\subset \overline{P^-_\epsilon}$ for $t\in[0,1]$.
\end{lemma}
{\bf Proof}: Since $K_c\backslash W=\emptyset$, we have $K_c\subset W$. By Lemma \ref{l5.1}, $I$
satisfies the (Ce)$_c$ condition, and so $K_c$ is compact. Then $2\delta:=dist(K_c,\partial W)>0$. For any $D\subset E$ and $a>0$, let $N_a(D):=\{u\in E: dist(u,D)<a\}$. Then $N_a(D)\subset W$. Moreover, since $I$ satisfies the (Ce)$_c$ condition, there exist $\epsilon_0, \alpha>0$ such that
$$(1+\|u\|)\|I'(u)\|\geq\alpha\ \text{for}\ u\in I^{-1}([c-\epsilon_0,c+\epsilon_0])\backslash {N_{\frac{\delta}{4}}(K_c)}.$$
By Lemmas \ref{l3.2} (2) and Lemma \ref{l3.5} (2), there exists $\beta>0$ such that
\begin{equation}\label{5.8}(1+\|u\|)\|u-B(u)\|\geq\beta\ \text{for}\ u\in I^{-1}([c-\epsilon_0,c+\epsilon_0])\backslash {N_{\frac{\delta}{4}}(K_c)}.\end{equation}
Without loss of generality, assume that $\epsilon_0<\frac{\beta\delta}{32(\sup_{u}K_c+2)}$. Let
$$h(u)=\frac{u-B(u)}{\|u-B(u)\|^2}\ \ \text{for}\ u\in E_0=E\backslash K,$$
and take a cut-off function $\zeta:E\rightarrow [0,1]$, which is locally Lipschitz continuous, such that
\begin{equation*}\aligned
\zeta(u)=\left\{ \begin{array}{lll}
0,\ & \text{if}\quad u\not\in I^{-1}[c-\epsilon',c+\epsilon'] \ \text{or}\ u\in N_{\frac{\delta}{4}}(K_c),\\
1, \ & \text{if}\quad u\in I^{-1}[c-\epsilon,c+\epsilon] \ \text{and}\ u\not\in N_{\frac{\delta}{2}}(K_c).
\end{array}\right.\endaligned
\end{equation*}
By Lemma \ref{l3.5}, $\zeta(\cdot)h(\cdot)$ is locally Lipschitz continuous on $E$.

Consider the following initial value problem
\begin{equation}\label{5.5}\aligned
\left\{ \begin{array}{lll}
\frac{d\tau}{dt}=-\zeta(\tau)h(\tau),\\
\tau(0,u)=u.
\end{array}\right.\endaligned
\end{equation}
For any $u\in E$,  by (\ref{5.8}) we have
\begin{equation*}
\|\zeta(u)h(u)\|\leq\frac{1}{\|u-B(u)\|}\leq\frac{1+\|u\|}{\beta}.\end{equation*}
 From Lemmas \ref{l5.2} and \ref{l5.3} it follows that the problem (\ref{5.5}) admits a unique solution $\tau(\cdot,u)\in C(\mathbb{R}^+, E)$.

 Define $\sigma(t,u)=\tau(16\epsilon t,u)$. Obviously, the conclusions (1) and (2) hold true. It suffices to check (3) and (4). To verify (3), let $u\in I^{c+\epsilon}\backslash W$. By Lemma \ref{l3.5}, $I(\tau(t,u))$ is decreasing in $t\geq0$. If there exists $t_0\in[0,16\epsilon]$ such that $I(\tau(t_0,u))<c-\epsilon$, then
$I(\sigma(1,u)=I(\tau(16\epsilon,u))<c-\epsilon$. Otherwise, for any $t\in[0,16\epsilon]$, there holds $I(\tau(16\epsilon,u))\geq c-\epsilon$. Then $\tau(t,u)\in I^{-1}[c-\epsilon,c+\epsilon]$ for any $t\in[0,16\epsilon]$. We claim that for any $t\in[0,16\epsilon]$, $\tau(t,u)\not\in N_{\frac\delta2}(K_c)$. Otherwise, there exists $t_1\in[0,16\epsilon]$ such that $\tau(t_1,u)\in N_{\frac\delta2}(K_c)$. Observe that $u\not\in N_{\delta}(K_c)$. Then
$$\text{there is} \ 0<t_2<t_1\ \text{such that}\ \tau(t_2,u)\in \partial(N_{\delta}(K_c)).$$ Hence
$$\aligned \frac{\delta}{2}&\leq\|\tau(t_2,u)-\tau(t_1,u)\|
\leq\int^{t_1}_{t_2}\|\tau'(s,u)\|ds\\
&\leq\int^{t_1}_{t_2}\frac{1}{\|\tau(s,u)-B(\tau(s,u))\|}ds\leq
\int^{t_1}_{t_2}\frac{1+\|\tau(s,u)\|}{\beta}ds\\
&\leq\frac{1+(\sup_{u}K_c+\delta)}{\beta}\int^{t_1}_{t_2}ds
\leq\frac{16\epsilon}{\beta}(\sup_{u}K_c+2),
\endaligned$$
which contradicts the fact that $\epsilon<\epsilon_0<\frac{\beta\delta}{32(\sup_{u}K_c+2)}$.
So $\zeta(\tau(t,u))\equiv1$ for $t\in [0,16\epsilon]$. From  Lemma \ref{l3.5} (2) and (3) it follows that
$$\aligned I(\tau(16\epsilon,u))&=I(u)-\int^{16\epsilon}_0\langle I'(\tau(s,u)),h(\tau(s,u))\rangle ds\\
&= I(u)-\int^{16\epsilon}_0\frac{\langle I'(\tau(s,u)), \tau(s,u)-B(\tau(s,u))\rangle}{\|\tau(s,u)-B(\tau(s,u))\|^2}ds\\
&\leq I(u)-\int^{16\epsilon}_0\frac18ds\leq c+\epsilon-2\epsilon=c-\epsilon.
\endaligned$$
Thus, the conclusion (3) yields. Finally, (4) is a consequence of  Lemma \ref{l3.5} (1), see \cite{LIUSUN} for a detailed proof.\ \ \ \ \ $\Box$

Arguing as Lemma \ref{l4.10} and using Lemma \ref{l3.1.0} (3) we have the following result.

\begin{lemma}\label{l3.9}If $\epsilon>0$ small enough, then $I(u)\geq\frac{\epsilon^2}{12}$ for $u\in \Sigma=\partial P^+_\epsilon\cap\partial P^-_\epsilon$, that is $c_*\geq\frac{\epsilon^2}{12}$.\end{lemma}

{\bf Proof of Theorem 1.2 (Existence part)} We will apply Theorem \ref{th1} to look for a sign-changing solution of (\ref{2.0.0}). By Lemma \ref{l5.2}, we know $\{P^+_\epsilon, P^-_\epsilon\}$ is an admissible family of invariant sets with respect
to $I$ at any level $c\in\mathbb{R}$. It suffices to verify assumptions (1)-(3) of
Theorem \ref{th1}. Let $v_1,v_2\in C^\infty_0(\mathbb{R}^3)\backslash\{0\}$ be such that $\text{supp}(v_1)\cap \text{supp}(v_2)=\emptyset$ and $v_1\leq0$, $v_2\geq0$. For $(t,s)\in \triangle$, let
$\varphi_0(t,s)=R(tv_1+sv_2),$
where the constant $R>0$ will be determined later. Obviously, for $t,s\in[0,1]$,
$\varphi_0(0,s)=Rsv_2\in P^+_\epsilon$ and $\varphi_0(t,0)=Rtv_1\in P^-_\epsilon$. Set $\rho=\min\{|tv_1+(1-t)v_2|_2:0\leq t\leq1\}>0$. Then $|u|_2\geq \rho R$ for $u\in \varphi_0(\partial_0\triangle)$. It follows from Lemma \ref{l4.9} that $\varphi_0(\partial_0\triangle)\cap M=\emptyset$. For any $u\in \varphi_0(\partial_0\triangle)$, we have
$$\aligned I(u)&=\frac{R^2t^2}2|\nabla v_1|^2_2+\frac12\int_{\mathbb{R}^3}V(x)f^2(Rtv_1)-\int_{\text{supp} v_1}G(f(Rtv_1))\\
&+\frac{R^2(1-t)^2}2|\nabla v_2|^2_2+\frac12\int_{\mathbb{R}^3}V(x)f^2(R(1-t)v_2)
-\int_{\text{supp}v_2}G(f(R(1-t)v_2))
.\endaligned$$
By ($g_4$) and Lemma \ref{l2.1} (5) we deduce $$\frac{f^2(s)}{s^2}\rightarrow0\ \text{and}\ \frac{G(f(s))}{s^2}=\frac{G(f(s))}{f^4(s)}\frac{f^4(s)}{s^2}\rightarrow+\infty,
\ \text{as}\ |s|\rightarrow+\infty.$$ we have $I(u)<0$ for large $R>0$. In view of Lemma \ref{l3.9} we have $c_*\geq\frac{\epsilon^2}{12}$ for small $\epsilon$. Hence
$\sup_{u\in \varphi_0(\partial_0\triangle)}I(u)<0<c_*.$
Applying Theorem \ref{th1} we know
(\ref{2.0.0}) admits a sign-changing solution $u_0\in E\backslash{(P^+_\epsilon \cup P^-_\epsilon)}$. Below arguing as the proof of Theorem 1.1 (Existence part), we show (\ref{2.0.0}) has a least energy sign-changing solution. \ \ \ \ $\Box$

\subsection{Multiplicity of sign-changing solutions}

 In order to apply Theorem \ref{th2}, we set $X=E$, $G=-id$, $J=I$ and $P=P^+_\epsilon$. Then $M=P^+_\epsilon\cap P^-_\epsilon$, $\Sigma=\partial P^+_\epsilon\cap \partial P^-_\epsilon$, and $W=P^+_\epsilon\cup P^-_\epsilon$. In this subsection, $g$ is assumed to be odd, and so $I$ is even. Since $K_c$ is compact, there exists a
symmetric open neighborhood $N$ of $K_c\backslash W$
such that $\gamma(\overline{N})<+\infty$. The counterpart of Lemma \ref{l4.11} under the (Ce) condition is as follows.

\begin{lemma}\label{l3.10} There exists $\epsilon_0>0$ such that, for $0<\epsilon<\epsilon'<\epsilon_0$, there
exists a continuous map $\sigma:[0,1]\times E\rightarrow E$ satisfying\\
\noindent(1) $\sigma(0,u)=u$ for $u\in E$;\\
\noindent(2) $\sigma(t,u)=u$ for $t\in[0,1]$, $u\not\in I^{-1}[c-\epsilon',c+\epsilon']$;\\
\noindent(3) $\sigma(1,I^{c+\epsilon}\backslash (N\cup W))\subset I^{c-\epsilon}$;\\
\noindent(4) $\sigma(t,\overline{P^+_\epsilon})\subset \overline{P^+_\epsilon}$ and $\sigma(t,\overline{P^-_\epsilon})\subset \overline{P^-_\epsilon}$ for $t\in[0,1]$;\\
\noindent(5) $\sigma(t,-u)=-\sigma(t,u)$ for $(t,u)\in[0,1]\times E$.\end{lemma}
{\bf Proof}: The proof of (1)-(4) are similar to those of Lemma \ref{l5.4}. Regarding (5), since  $I$ is even and the operator $B$ is odd in Lemma \ref{l3.5}, we infer
$\sigma$ is odd in $u$.\ \ \ \ $\Box$

{\bf Proof of Theorem 1.2} ({\bf Multiplicity part}) In view of Lemma \ref{l3.10} we know
$P^+_\epsilon$ is a $G-$admissible invariant set with respect to $I$ at
level $c\in\mathbb{R}$. To apply Theorem \ref{th2}, we are now constructing $\phi_n$. For any $n\in\mathbb{N}$, choose $\{v_i\}^n_1\subset C^\infty_0(\mathbb{R}^3)\backslash\{0\}$ such that $supp(v_i)\cap supp(v_j)=\emptyset$ for $i\neq j$. Define $\varphi_n\in C(B_n, E)$
as
$$\varphi_n(t)=R_n\Sigma^{n}_{i=1}t_iv_i,\quad\ t=(t_1,t_2,...,t_n)\in B_n,$$
where $R_n>0$. For $R_n$ large enough, it is easy to check that all the assumptions of Theorem \ref{th2} are satisfied. Therefore,  (\ref{2.0.0}) admits infinitely many sign-changing solutions.\ \ \ \ $\Box$

\section{Proof of Theorem 1.3}
\renewcommand{\theequation}{5.\arabic{equation}}
In this section, we consider equation (\ref{1.0.1}). By the change $f$, (\ref{1.0.1}) is changed into
\begin{equation}\label{6.1}-\Delta u+V(x)f(u)f'(u)=f^3(u)f'(u).\end{equation}
 In this case, the associated functional $I$ merely satisfies the (Ce) condition. So we use the deformation lemmas as in Sections 4.3 and 4.4 to show the existence and multiplicity of sign-changing solutions of (\ref{6.1}). However, in this case, the condition
\begin{equation}\label{6.2}\lim_{|t|\rightarrow+\infty}\frac{G(t)}{t^4}=+\infty,\end{equation}
 is not satisfied, then it is necessary to make  suitable modifications in the arguments of Lemma \ref{l5.1} and Theorem 1.2. And the other results and their proof are similar to those in Section 4, so they will not be written out.

\begin{lemma}\label{l6.1} For any $c\in\mathbb{R}$, $I$ satisfies the (Ce)$_c$ condition.\end{lemma}
{\bf Proof}: Assume that $\{u_n\}\subset E$ satisfies
$I(u_n)\rightarrow c$ {and} $(1+\|u_n\|)\|I'(u_n)\|\rightarrow0.$
 Let $w_n=\frac{f(u_n)}{f'(u_n)}$. By (\ref{3.0.1}) we deduce
$$ \aligned c+o_n(1)&=I(u_n)-\frac14\langle I'(u_n),w_n\rangle=\frac{1}{4}\int_{\mathbb{R}^3}(f'(u_n))^2|\nabla u_n|^2+\frac14\int_{\mathbb{R}^3} V(x)f^2(u_n)\\&=\frac{1}{4}\int_{\mathbb{R}^3}|\nabla f(u_n)|^2+\frac14\int_{\mathbb{R}^3} V(x)f^2(u_n).\endaligned$$
Therefore, $\{f(u_n)\}$ is bounded in $E$ and then bounded in $L^s(\mathbb{R}^3)$ with $2\leq s\leq6$. Note that
 $$c+o_n(1)=I(u_n)=\frac12|\nabla u_n|^2_2+\frac12\int_{\mathbb{R}^3}V(x)f^2(u_n)-\frac14|f(u_n)|^4_4.$$
So $\{|\nabla u_n|^2_2\}$ is bounded and $\{u_n\}$ is bounded in $E$ using (\ref{2.3.4}).  In the same way as (\ref{4.2.4}) we have $u_n\rightarrow u$ in $E$. Thus $I$ satisfies the (Ce)$_c$ condition.\ \ \ \ $\Box$

{\bf Proof of Theorem 1.3} We take similar arguments of Theorem 1.2, and without the condition (\ref{6.2}), we make some modification as follows. In applying  Theorem \ref{th1}, we need to choose $v_1,v_2\in C^\infty_0(\mathbb{R}^3)\backslash\{0\}$ to be such that  $$|\nabla v_i|^2_2<|v_i|^2_2,\ i=1,2,\ \text{supp}(v_1)\cap \text{supp}(v_2)=\emptyset,\ \text{and} \ v_1\leq0, v_2\geq0.$$ For $(t,s)\in \triangle$, let
$\varphi_0(t,s)=R(tv_1+sv_2),$
where $R>0$ will be determined later. For any $u\in \varphi_0(\partial_0\triangle)$, we have
$$\aligned I(u)&=\frac{R^2t^2}2|\nabla v_1|^2_2+\frac12\int_{\mathbb{R}^3}V(x)f^2(Rtv_1)-\frac14\int_{\mathbb{R}^3}|f(Rtv_1)|^4\\
&+\frac{R^2(1-t)^2}2|\nabla v_2|^2_2+\frac12\int_{\mathbb{R}^3}V(x)f^2(R(1-t)v_2)
-\frac14\int_{\mathbb{R}^3}|f(R(1-t)v_2)|^4
.\endaligned$$
In view of  Lemma \ref{l2.1} (5) we deduce $$\frac{f^2(s)}{s^2}\rightarrow0\ \text{and}\ \frac{|f(s)|^4}{s^2}\rightarrow2,
\ \text{as}\ |s|\rightarrow+\infty.$$ Combining with the fact that $|\nabla v_i|_2<|v_i|_2$, we have $I(u)<0$ for large $R>0$. Then arguing as Theorem 1.2 (existence part), we infer that (\ref{6.1}) has a least energy sign-changing solution.

In applying Theorem \ref{th2}, as above we need to  construct a different $\phi_n$ from that in the proof of Theorem 1.2 (multiplicity part). For any $n\in\mathbb{N}$, let $\{v_i\}^n_1\subset C^\infty_0(\mathbb{R}^3)\backslash\{0\}$ be such that $$|\nabla v_i|^2_2<|v_i|^2_2,\ \text{for all}\ i=1,2,...,n,\ \text{and}\ \text{supp}(v_i)\cap \text{supp}(v_j)=\emptyset\ \text{for}\ i\neq j.$$ Define $\varphi_n\in C(B_n, E)$
as
$\varphi_n(t)=R_n\Sigma^{n}_{i=1}t_iv_i$, where $t=(t_1,t_2,...,t_n)\in B_n$
and $R_n>0$. For $R_n$ large enough, as above one easily check that all the assumptions of Theorem \ref{th2} are satisfied. Hence, equation (\ref{6.1}) has infinitely many sign-changing solutions.\ \ \ $\Box$


\begin{thebibliography}{100}

\bibitem{BLW} T. Bartsch, Z.L. Liu, T. Weth, Nodal solutions of a $p$-Laplacian equation, {\it Proc. Lond. Math. Soc.}
{\bf 91} (2005), 129-152.

\bibitem{BERE} H. Berestycki, P.L. Lions,  Nonlinear scalar field equation I.
Existence of a ground state,  {\it Arch. Rational Mech. Anal.} {\bf 82} (1983), 313-345.
\bibitem{BL} L. Br\"{u}ll, H. Lange, Solitary waves for quasilinear Schr\"{o}dinger equations, {\it Exposition Math.} {\bf 4} (1986), 279-288.

\bibitem{ChenS}  X.L. Chen, R.N. Sudan, Necessary and sufficient conditions for self-focusing of short ultraintense laser pulse in underdense
plasma, {\it Phys. Rev. Lett.}  {\bf70} (1993), 2082-2085.

\bibitem{CJ} M. Colin, L. Jeanjean, Solutions for a quasilinear Schr\"{o}dinger equation: A dual approach, {\it Nonlinear Anal.} {\bf56} (2004), 213-226.

\bibitem{DPW1} Y.B.  Deng, S.J. Peng, J.X. Wang, Infinitely many sign-changing solutions for quasilinear Schr\"{o}dinger equations in $\mathbb{R}^N$, {\it Commun. Math. Sci.} {\bf9} (2011), 859-878.

 \bibitem{DPW2} Y.B.  Deng, S.J. Peng, J.X. Wang, Nodal soliton solutions for quasilinear Schr\"{o}dinger equations with critical exponent, {\it J. Math. Phys.} {\bf54} (2013), 011504, 27 pp.

\bibitem{DS1} J.M. do \'{O}, U. Severo, Solitary waves for a class of quasilinear Schr\"{o}dinger equations in
dimension two, {\it Calc. Var. Partial Differential Equations} {\bf38} (2010), 275-315.

\bibitem{Jeanjean}  L. Jeanjean, On the existence of bounded palais-smale sequences and application to a landesman-lazer-type problem set on $\mathbb{R}^N$, {\it Proc. Roy. Soc. Edinburgh Sect. A }{\bf129} (1999), 787-809.

\bibitem{Ku} S. Kurihara,  Large-amplitude quasi-solitons in superfluid films, {\it J. Phys. Soc. Jpn.} {\bf50} (1981), 3262-3267.


 \bibitem{LPT}  H. Lange, M. Poppenberg, H. Teismann, Nash-Moser methods for the solution of quasilinear Schr\"{o}dinger equations,
{\it Comm. Partial Differential Equations} {\bf24} (1999), 1399-1418.

\bibitem{Ligongbao} G.B. Li, C.H. Wang, The existence of a nontrivial solution to a nonlinear elliptic problem of linking type without the Ambrosetti-Rabinowitz condition, {\it Ann. Acad. Sci. Fenn-m} {\bf36} (2011), 461-480.

\bibitem{LLW} J.-Q. Liu, X.Q. Liu, Z.-Q. Wang, Multiple mixed states of nodal solutions for nonlinear Schr\"{o}dinger
systems, {\it Calc. Var. Partial Differential Equations} {\bf52} (2015), 565-586.

\bibitem{LLW2} J.-Q. Liu, X.Q. Liu, Z.-Q. Wang, Multiple sign-changing solutions for quasilinear elliptic equations via perturbation method, {\it Comm. Partial Differential Equations}, {\bf 39}(2014), 2216-2239.


\bibitem{LWW} J.-Q. Liu, Y.-Q. Wang, Z.-Q. Wang, Solutions for quasilinear Schr\"{o}dinger equations via the Nehari method, {\it Comm. Partial Differential Equations} {\bf29} (2004), 879-901.

 \bibitem{LIUSUN} Z.L. Liu, J.X. Sun, Invariant sets of descending flow in critical point theory with applications to nonlinear
differential equations, {\it J. Differential Equations} {\bf172} (2001), 257-299.

 \bibitem{LWZ} Z.L. Liu, Z.-Q. Wang, J.J. Zhang, Infinitely many sign-changing solutions for the nonlinear
Schr\"{o}dinger-Poisson system, {\it Ann. Mat. Pur. Appl.} {\bf195} (2016), 775-794.

\bibitem{LZ} Z.S. Liu, Y.J. Lou, J.J. Zhang, A perturbation approach to studying sign-changing solutions of Kirchhoff equations with a general nonlinearity, arXiv:1812.09240.

\bibitem{LS} Z.S. Liu, G. Siciliano, A perturbation approach for the Schr\"odinger-Born-Infeld system: solutions in the subcritical and critical case, {\it J. Math. Anal. Appl.} {\bf503} (2021), 125326, 22 pp.


\bibitem{PS} M. Poppenberg, K. Schmitt, Z.-Q. Wang, On the existence of soliton solutions to
quasilinear Schr\"{o}dinger equations, {\it Calc. Var. Partial Differential Equations} {\bf14} (2002), 329-344.


\bibitem{Rap} P. Rabinowitz, On a class of nonlinear Schr\"{o}dinger equations, {\it Z. Angew. Math. Phys.}  {\bf43} (1992), 270-291.


\bibitem{WM1} M. Willem, Minimax Theorems,
Progr. Nonlinear Differential Equations Appl, vol.24,
Birkh$\ddot{a}$user, Basel, 1996.

\bibitem{XC} L.P. Xu, H.B. Chen, Ground state solutions for quasilinear Schr\"{o}dinger equations via Pohozaev manifold in Orlicz space, {\it J. Differential Equations} {\bf265} (2018), 4417-4441.

 \bibitem{YANG1} M.B. Yang, C.A. Santos, J.Z. Zhou,
Least energy nodal solutions for a defocusing Schr\"{o}dinger equation with supercritical exponent,
{\it Proc. Edinb. Math. Soc.} {\bf62} (2019), 1-23.




\bibitem{ZL} W. Zhang, X.Q. Liu, Infinitely many sign-changing solutions for a quasilinear elliptic equation in $\mathbb{R}^N$, {\it J. Math. Anal. Appl.} {\bf427} (2015), 722-740.
\bibitem{Zhongchengkui} C.K. Zhong, X.L. Fan, W.Y. Cheng, An introduction to nonlinear functional analysis, Lanzhou Univ. Press, Lanzhou, 1998.





\end{thebibliography}
\end{document}